\documentclass[12pt]{article}
 
\usepackage[bottom=30mm, top=30mm, left=25mm,right=25mm]{geometry}
\usepackage{setspace}
\usepackage{cite}
\usepackage{amssymb,amsmath}
\usepackage{bbm}
\usepackage{pdflscape}
\usepackage{subfig}
\usepackage{graphicx}
\usepackage{xcolor}

\newcommand{\cb}[1]{{\color{blue} #1}}
\newcommand{\crr}[1]{{\color{red} #1}}

\newcommand{\bR}{\mathbb{R}}

\usepackage{hyperref}

\begin{document}
\numberwithin{equation}{section}
\setstretch{1.2}
\setlength{\parskip}{1ex plus 0.5ex minus 0.2ex}
\setlength{\parindent}{0pt}

~\vspace{1cm}~
\begin{center}
{\Large\textbf{Learning Topological Invariance}}\\[10mm]

James Halverson$^{a,c,}$\footnote{\href{mailto:j.halverson@northeastern.edu}{j.halverson@northeastern.edu}} and Fabian Ruehle$^{a,b,c,}$\footnote{\href{mailto:f.ruehle@northeastern.edu}{f.ruehle@northeastern.edu}}\\[10mm]
{
	{\it ${}^{\text{a}}$ Department of Physics, Northeastern University, Boston, MA 02115, USA}\\[.5em]
	{\it ${}^{\text{b}}$ Department of Mathematics, Northeastern University, Boston, MA 02115, USA}\\[.5em]
	{\it ${}^{\text{c}}$ NSF Institute for Artificial Intelligence and Fundamental Interactions}\\[.5em]
}
\end{center}
\setcounter{footnote}{0}
\vspace{24pt}
\begin{abstract}
\noindent Two geometric spaces are in the same topological class if they are related by certain geometric deformations. We propose machine learning methods that automate learning of topological invariance and apply it in the context of knot theory, where two knots are equivalent if they are related by ambient space isotopy. Specifically, given only the knot and no information about its topological invariants, we employ contrastive and generative machine learning techniques to map different representatives of the same knot class to the same point in an embedding vector space. An auto-regressive decoder Transformer network can then generate new representatives from the same knot class. We also describe a student-teacher setup that we use to interpret which known knot invariants are learned by the neural networks to compute the embeddings, and observe a strong correlation with the Goeritz matrix in all setups that we tested. We also develop an approach to resolving the Jones Unknot Conjecture by exploring the vicinity of the embedding space of the Jones polynomial near the locus where the unknots cluster, which we use to generate braid words with simple Jones polynomials.
\end{abstract}
\clearpage
\tableofcontents
\clearpage

\section{Introduction}

In topology, spaces are said to be equivalent if they can be continuously deformed into one another. Verifying that this is the case for any two particular spaces is in general a difficult problem, but falsifying it is not: if the spaces differ on a single topological invariant, they are not topologically equivalent. In recent years, machine learning techniques have been utilized for both verifying and falsifying topological equivalence. For instance, in knot theory studies where topological invariants are learned with supervised learning include \cite{hughes,Jejjala:2019kio,Craven:2021ckk,Davies2021,Craven:2022cxe, lindsay2025learnabilityknotinvariantsrepresentation}. On the other hand, rigorous verification of topological equivalence \cite{gukov2020learningunknot, applebaum2024unknottingnumberhardunknot} to the unknot and concordance \cite{Gukov:2023kvx,shehper2025makesmathproblemshard} to the unknot was learned with reinforcement learning (RL); the latter ruled out many proposed counterexamples to the smooth Poincar\'e conjecture in four dimensions. Both the equivalence and concordance RL games rely on theorems that establish manipulations of knots sufficient to demonstrate equivalence or concordance, e.g., two knots are equivalent iff their projections can be connected by a sequence of Reidemeister moves~\cite{Reidemeister1927}. The sequences may be long, however, leading to the difficulty of the game. 

In this paper, we probe whether ML can learn topological invariance given many representatives of different equivalence classes, without any a priori knowledge of theorems related to topological invariants that distinguishes these classes. Moreover, we are interested in interpreting the result, i.e., extracting any known invariants that the Neural Network (NN) implicitly learned in order to exhibit topological invariance. If successful, the trained network should cluster topologically equivalent objects in the embedding dimension and be able to generate new representatives of the same knot class. 

We utilize these techniques rather than direct supervised classification because the latter might learn equivariance rather than invariance. In the equivariant case, there could be multiple clusters or entire submanifolds in the embedding space that belong to the same equivalence class, which are transformed into one another under actions that do not change the topological class. These clusters or submanifolds would then have to be arranged in the embedding space such that they ``project'' to the same output after applying the final layer. Moreover, classification requires fixing the number of equivalence classes at design time.

We study topological invariance learning in the context of knot theory. We use braid words to represent the knots. A braid word in the non-Abelian braid group given by a sequence of Artin braid generators, which may be represented as a string of integers whose absolute value is bounded by the number of strands in the braid minus one. We explore different losses (contrastive and generative learning) as well as different NN architectures. For contrastive learning, we use a simple feed-forward NN and a convolutional NN (CNN). The convolutional layer is used to obtain invariance under one of braid group symmetries. For the losses, we employ two different online contrastive learning losses: semi-hard triplet loss and a centroid loss that we develop specifically for this problem. For the generative approach, we train an encooder-decoder Transformer~\cite{NIPS2017_3f5ee243} pair to generate braid words. In more detail, the Transformers are trained on many different knot classes with many representatives each, requiring that each representative of a given class map to a fixed representative at output. We test whether the decoder constructs knots that are topologically equivalent to the input knot by comparing their hyperbolic knot volumes. Remarkably, both methods are able to learn meaningful embeddings of knot classes that encode topological invariance. 

We are also interested in extracting the invariants from the embedding dimension that are learned by the NN. To do this, we use a student-teacher setup. The teachers are large NNs that learn invariance from contrastive or generative learning. We use their learned braid embeddings as labels for simpler student networks, which are tasked to reproduce these labels from known knot invariants. We take the degree (more precisely, the MSE loss) to which the student NNs can learn a map from the known knot invariants to the teacher embeddings as a measure for ``how much'' of a given knot invariant is used by the teacher to compute the embedding. We observe that across all different architectures and methods (contrastive and generative), the NNs seem to learn invariants based on the Goeritz matrix.

Finally, we use the generative learning approach to search for counter-examples to the Jones Unknot Conjecture, which states that a knot is the unknot if and only if it has a trivial Jones polynomial. We train a Transformer to map Jones polynomials to braid word representatives. By searching the vicinity of knots with simple Jones polynomials in the embedding space of the Transformer, we aim to generate other knots with simple Jones polynomials that are distinct from the unknot.

This paper is organized as follows. In Section \ref{sec:ideas} we describe the conceptual learning setup, including the learning of topological invariance, extracting correlations with known topological invariants, and the Jones unknot problem. In Section~\ref{sec:expts} we present the experimental results, including data generation, learning invariance, extracting invariants, and an analysis of the Jones Unknot Conjecture. We conclude in Section~\ref{sec:conclusions}. Appendix~\ref{app:GeneratedLinks} illustrates the different type of knot distributions we get from the  generative approaches for knots with simple Jones polynomials, and Appendix~\ref{app:BeamSearch} illustrates in a few examples how different knots that generated by the Transformer decoder are  related by Markov moves and braid relations.

Note added: while finalizing this manuscript, the preprint \cite{braghetto2025variationalautoencodersunderstandknot} appeared, which has some overlapping ideas with this work. Also, Gujral and Hughes are currently finishing a preprint~\cite{Gujral:2025aaa} that has some overlap with this work.

\section{Conceptual Learning Setup\label{sec:ideas}}

\subsection{Learning Invariance\label{sec:idea_invariance}}
In this section we present an algorithm that learns equivalence relations.  More precisely, given a set $S$ and an equivalence relation on the elements of the set, we partition the set into equivalence classes and train a neural network to map elements of the same equivalence class to the same point in an $n$-dimensional embedding space $\mathbbm{R}^n$. By learning to map elements from the same equivalence class to the same point, the NN is learning invariance under the equivalence relation. By mapping elements from distinct equivalence classes to distinct points, we ensure that the NN is learning invariance under the equivalence relation instead of a stronger notion of equivalence. In particular, this ensures that the NN is not learning a trivial solution where every element is mapped to the same point in $\mathbbm{R}^n$. 

In the following, we illustrate this idea by learning topological invariance in knot theory, but the idea is applicable more broadly.

A knot is an embedded $S^1$ in $S^3$ up to ambient isotopy. Hence, in our setup, the set $S$ is the set of all embedded circles into a three-sphere, and two knots are equivalent if and only if one can be continuously deformed into the other without breaking it. Schematically, let us denote the set of all knots by $\mathcal{K}$, the different equivalence classes of knots by $\mathcal{K}^i$, and the equivalent knots in each class by $K^i_\alpha\in\mathcal{K}^i$. Evidently, we have an infinite set of equivalence classes $\mathcal{K}^i$, each of which has an infinite set of representatives. We set up a NN $f_\theta$, which is a parameterized map from (a representation of) $\mathcal{K}$ into $\mathbbm{R}^n$ such that
\begin{align}
\label{eq:KnotEquivalence}
\begin{split}
f_\theta:~~	&\mathcal{K}^i_{\phantom{\alpha}}\to\mathbbm{R}^n\,,\\
			&K^i_\alpha \mapsto m_i\in\mathbbm{R}^n\qquad\forall\alpha\,, \\
            &m_i \neq m_j \qquad\qquad\quad \forall i \neq j,
\end{split}
\end{align}
i.e., all knots from the same equivalence class are mapped to the same embedding vector. The $m_i$ are then knot invariants. 

We follow two approaches to train the NN. The first is using contrastive learning, where we engineer a loss function that ensures relation~\eqref{eq:KnotEquivalence} directly on the embedding space. The second trains a generative model that instead maps
\begin{align}
\label{eq:KnotGeneration}
\begin{split}
f_\theta:~~	&\mathcal{K}^i_{\phantom{\alpha}}\to\mathcal{K}^i\,,\\
			&K^i_\alpha \mapsto K^i_1 \qquad\forall\alpha\,.
\end{split}
\end{align}
Here, we choose a (random) representative $K^i_1$ from each equivalence class $\mathcal{K}^i$, and train the neural network to reproduce this representative as output from any input representative of the same class (we will see how we embed $K^i_\alpha$ in $\mathbbm{R}^n$ in the subsection below). Again, the NN learns invariance under equivalence of inputs. We call this technique generative learning because we utilize a Transformer to generate the output braid representative, but we note that it could also be used with a simpler neural network such as a CNN or an MLP.

In the generative approach, the last layer of the NN before the output layer can be used as the embedding layer and can be thought of as the space $\bR^{n}$ in \eqref{eq:KnotEquivalence}. Of course, both approaches could also be combined, using a generative network for which we also impose a contrastive learning loss in the embedding layer.

\subsubsection{Knot representations}
\label{sec:KnotReps}
In order to input a knot to a NN, we need to represent it as a vector in $\mathbbm{R}^{n_\text{in}}$. The above description would be most directly implemented by sampling $P$ points along the embedded $S^1$ that describes the knot and representing it as a polygonal curve in $\mathbbm{R}^{3P}$. Equivalent representations of the knot are obtained by deforming the embedded $S^1$. We did not attempt this here but plan to return to it in the future. Instead, we represent the knot as a braid. 

The $n$-stranded braid group $B_n$ is a non-Abelian group with $n-1$ generators $\sigma_i$ and their inverses $\sigma_i^{-1}$, $i=1,\ldots,n-1$. The generator $\sigma_i$ moves strand $i$ over strand $i+1$, and the inverse moves strand $i$ under strand $i+1$. It is easy to check that braid actions are associative, closed under composition, and have a neutral element where strands are not moved, i.e., they form a group. The group generators satisfy the commutation relation
\begin{align}
\sigma_i\circ\sigma_j=\sigma_j\circ\sigma_i\,,\qquad |i-j| > 1\,,
\end{align}
and the relation
\begin{align}
\sigma_i\circ\sigma_{i+1}\circ\sigma_{i}=\sigma_{i+1}\circ\sigma_{i}\circ\sigma_{i+1}\,.
\end{align}

Alexander showed that every knot can be written as a closed braid~\cite{Alexander:1923aaa}, which means that the end points of the strands are identified. This identification means that generators on one end of the braid can be moved to the other end. More formally, two braid words $w$ that differ by conjugation, $w\to g w g^{-1}$ for some $g\in B_n$, describe the same knot. Moreover, adding a new generator $\sigma_{n}^{\pm1}$ to a braid word $w\in B_n$ gives an element $w\to w\sigma_{n}^{\pm1} \in B_{n+1}$ that still describes the same knot (at the level of the knot, this ``stabilization move'' introduces an additional loop somewhere in the knot that can just be straightened out again). These two operations are called Markov moves, and two knots are related by ambient isotopy iff one braid representation can be brought to the other using only the two Markov moves and the two braid relations.

We will use the braid representation to encode the knot such that it can be input into the NN. A straight-forward way, employed e.g.\ in~\cite{gukov2020learningunknot}, is to encode a knot  represented by a braid word $w$ as a vector in $\mathbbm{Z}^{|w|}$, where we map $\sigma_i$ to $i$ and $\sigma_i^{-1}$ to $-i$. For large enough $n$, one can additionally scale the input to have zero mean (which it should roughly have anyways, since neither $\sigma$ nor $\sigma^{-1}$ are preferred) and unit variance.

Another way to embed the knot is to tokenize the input. In this procedure, we first one-hot-encode each braid generator of $B_n$, i.e., we assign them to a $2(n-1)$-dimensional vector where
\begin{align}
\begin{split}
\sigma_1\mapsto(1,0,0,0,\ldots,0,0)\,,&\qquad \sigma_1^{-1}\mapsto(0,1,0,0,\ldots,0,0)\,,\\
\sigma_2\mapsto(0,0,1,0,\ldots,0,0)\,,&\qquad\ldots\,,\\
\sigma_n\mapsto(0,0,0,0,\ldots,1,0)\,,&\qquad\sigma_n^{-1}\mapsto(0,0,0,0,\ldots,0,1)\,.
\end{split}
\end{align}
In the next step, this is converted to a vector in $\mathbbm{R}^{n_\text{emb}}$ by multiplying each one-hot-encoded vector by an $n_\text{emb}\times 2n$ matrix. 

\subsubsection{Neural network architectures}
Besides the two different braid encodings discussed above, we explored different network architectures for both the contrastive and the generative learning objective. 

\textbf{Autoencoders.} We tried (variational) autoencoders~\cite{kingma2013auto} with standard feed-forward neural networks or Kolmogorov-Arnold Networks~\cite{Liu:2024swq} for the generative task, but the results were subpar compared to Transformers, so we will focus on the latter. We also used just the encoder part (with an MLP or CNN architecture) with a contrastive loss in the bottleneck layer. In that case, we perform a Principal Component Analysis (PCA) in the bottleneck to compute the effective embedding dimension. For the contrastive loss, we find that essentially all of the explained variance in the PCA gets concentrated in 2-3 dimensions, as well will further discuss in Section~\ref{sec:expts}.

\textbf{CNNs.} We used CNNs to enforce invariance of the neural net under closure of the braid word (i.e., under conjugation by the generators at the end of the braid word). In more detail, we set up the first layer to be a convolutional 1D layer with filter size the length of the braid word, stride 1, and periodic boundary conditions. We denote the braid word by $w_B=w_1 w_2\ldots w_n $ where each $w_i$ is one of the braid generators or their inverse, the weights of the convolutional filter by $\theta_j$, $j=1,\ldots,n$, and the result of the convolution as $\vec{v}\in\mathbbm{R}^n$. We then get
\begin{align}
\begin{split}
v_1 &= \theta_1 w_1 + \theta_2 w_{2\phantom{+j}} + \ldots + \theta_{n-1} w_{n-1} ~\,+ \theta_n w_n\,,\\
v_2 &= \theta_1 w_2 + \theta_2 w_{3\phantom{+j}} + \ldots + \theta_{n-1} w_{n\phantom{-1}} ~\,+ \theta_n w_1\,,\ldots\\
v_j &= \theta_1 w_{j} + \theta_2 w_{j+1} + 	    \ldots + \theta_{(n-j+1)} w_{n} + \theta_{n-j+2} w_{1} + \ldots\,.
\end{split}
\end{align}
Hence, the resulting vector $\vec{v}$ contains a sequence of all possible conjugations by $w_1 w_2 \ldots w_j$, $j=1,\ldots n$.

\textbf{Transformers.} The Transformer consists of an encoder and a decoder that together map sequences to sequences. We use so-called ``tokenized" braid words as inputs and outputs, and introduce a special token that signals the end of the braid word. The latter is needed since it is an auto-regressive model, which means that the decoder is presented with its own output plus information about the input from the encoder to predict the token that corresponds to the next braid word generator, or the end of sequence token to indicate that the braid word is finished. The Transformer architecture itself consists of an encoder and a decoder, which are built from stacks of so-called attention heads and fully connected layers, see~\cite{Vaswani:2017Attention} for details. 

The auto-regressive generation process works as follows. The encoder receives tokenized braid words. It embeds the one-hot-encoded generators as described above and adds a positional embedding that indicates at which position in the braid word the generator is located. Fully connected layers are agnostic of the order of the input features, but the order clearly matters in natural language, or in our case for non-Abelian groups such as the braid group. Of course, the NN can in principle learn that the input is ordered, but the temporal embedding is intended to provide an additional signal which the network can pick up that would tell it the ordering. These embedded tokens are then fed through the multi-head attention encoder to produce vectors in $\mathbbm{R}^{n_\text{enc}}$. For the decoder, we start with the empty string (i.e., the end of sequence token), but also provide the vectors generated by the encoder. The final layer of the decoder is a softmax layer in $\mathbbm{R}^{2n+1}$. The $2n+1$ entries represent the probability with which one of the $n$ braid generators, their inverse, or the end of sequence indicator should be added to the output braid word. We select the token with the highest probability and add it to the generated word. After this step, the decoder is presented with the braid word it has generated so far, plus the input from the encoder, to predict the next generator (or end of sequence token) to be added to the braid word. This procedure is repeated until the decoder predicts the end of the sequence or a maximal length is reached. 

Instead of greedily selecting the most probable next token, we can also perform a so-called beam search where we add the $k$ most probable next tokens to the output, or where we sample from the probability distribution produced by the decoder output. For example,  a beam search can provide multiple different braid representatives of a fixed topological class.

\subsubsection{Contrastive learning}
A popular supervised contrastive learning~\cite{Chopra:2005Contrastive} loss is the so-called triplet loss~\cite{Schroff:2015Facenet} and its variations. The idea of triplet loss is that the NN is trained with a triple $(A,P,N)$ where $A$ is an anchor, $P$ is a positive sample and $N$ is a negative sample. The anchor and the positive are taken to be elements from the same equivalence class $\mathcal{K}^i$, while the negative is a sample from a different equivalence class $\mathcal{K}^j$ with $j \neq i$.

Triplet loss then minimizes the (Euclidean) distance between the anchor and the positive while maximizing the (Euclidean) distance between the anchor and the negative,
\begin{align}
\mathcal{L}(A,P,N)=\max(||f_\theta(A)-f_\theta(P)||^2-||f_\theta(A)-f_\theta(N)||^2+\kappa,0)\,.
\end{align} 
The first term pulls the positive towards the anchor, while the second term pushes the negative away from the anchor. Note that with just the first two terms, the loss would not be positive, which is why the maximum is included. With the maximum, the loss will be zero if the distance from the anchor to the positive is less than the distance from the anchor to the negative. To enforce better clustering, a constant $\kappa\in\mathbbm{R}^+$ is added to the loss, such that now the distance between the anchor and the negative has to be at least $\kappa$ plus the distance from the anchor and the positive to result in zero loss. The best results are obtained from semi-hard triplets, which are triplets where the negative is farther from the anchor than the positive (but by less than $\kappa$, such that the loss is still non-zero). This requires online learning, where semi-hard triples are chosen on the fly and only those are used during training. 

The above loss requires the choice of (random) anchors, which we found theoretically unpleasing. Also, elements in the embedding space can be mapped to different points (within the same decision region) as long as they are at least a distance $\kappa$ from the decision boundary. However, we demanded in~\eqref{eq:KnotEquivalence} that all elements get mapped to the same $m_i\in\mathbbm{R}^n$. 

To circumvent these points, we introduce a centroid loss. The idea is to replace the anchor by the centroid computed from all instances from the same equivalence class. This also requires online learning, where the centroids are updated after each gradient step. In more detail, let
\begin{align}
c_i = \frac1N \sum_{\alpha=1}^N f_\theta(K^i_\alpha)
\end{align}
be the centroids of the equivalence classes $\mathcal{K}^i$ in $\mathbbm{R}^n$. These will become the knot invariants $m_i$. If the number of different equivalence classes is $M$, and the number of representatives of each equivalence class is $N$ as above, we define the centroid loss as
\begin{align}
\mathcal{L}(c_i,K^j_\alpha)=\frac1M\sum_{i=1}^M\frac1N\sum_{\alpha=1}^N||c_i-f_\theta(K^i_\alpha)||^2\,.
\end{align}
Note that this loss is zero iff all elements from the same equivalence class get mapped to the same point $c_i$. This simple version does not enforce a minimal distance between centers, which can be enforced by adding a repulsion term for cluster centers,
\begin{align}
\mathcal{L}(c_i,K^j_\alpha)=\frac1M\sum_{i=1}^M\frac1N\sum_{\alpha=1}^N||c_i-f_\theta(K^i_\alpha)||^2 + \lambda \sum_{i,j=1}^M \max(\kappa-||c_i-c_j||,0)\,.
\end{align}
The second term will be zero if all cluster centers are at least a distance of $\kappa\in\mathbbm{R}^+$ apart, and $\lambda\in\mathbbm{R}^+$ is a relative weight between the attraction term towards cluster centers and the repulsion term of clusters of different equivalence classes.

\subsubsection{Generative learning}
We obtained the best results from Transformers with tokenized sequences, standard trigonometric positional embedding, and cross-entropy loss. As explained above, we pick an arbitrary representative from each equivalence class as the target for the decoder. However, this presents a subtlety for evaluating performance of the Transformer, since we do not know which representative it will generate on the test or validation set. Hence, we cannot label the unseen data. To decide whether the Transformer produced a correct representative, we compute knot invariants of the input and the output knot and compare them. We restrict ourselves to comparing the hyperbolic knot volume, which is a simple scalar invariant that can be computed efficiently using the python package snappy~\cite{SnapPy}. While this is only defined for hyperbolic knots, and there exist inequivalent knots with the same volume, most knots (at least in the range of numbers of crossings we are considering) are hyperbolic and volume collisions are rare; i.e., volumes are sufficient to distinguish most knot classes. In fact, it is highly unlikely that the Transformer would output a knot from a different equivalence class that happens to have the same volume as the knots from the input equivalence class. 

\subsection{Extracting Invariants\label{sec:idea_invariants}}
After successful training, the NN will have learned to embed all equivalent knots to the same point in the embedding dimension, i.e., the embedding vector $m_i\in\mathbbm{R}^n$ encodes a knot invariant. Note that the knot invariant is not necessarily $n$-dimensional though, since it is conceivable that all knots live on a (high-codimension) submanifold inside $\mathbbm{R}^n$. Indeed, for the contrastive learning tasks we found that the NN effectively only uses 2 or 3 dimensions, irrespective how big we make the embedding dimension. We find the intrinsic dimension of the embedding space by performing a principal component analysis (PCA) and cutting off all principal components (PCs) above a certain threshold for the total explained variance. 

Once we have the PCs that encode the knot invariants learned by the NN, we want to extract which (combination of) known knot invariants the NN learned, and potentially uncover new invariants uncovered by the NN. In order to compare against known knot invariants, we employ a student-teacher-like setup. Using contrastive or generative methods to train the teacher, we compute for each equivalence class $\mathcal{K}_i$ the PCs of $m_i$ up to an explained variance threshold. We then train a (smaller) student NN in a supervised fashion, where the inputs are known knot invariants and the outputs are the PCs from the teacher. We look at several invariants, including
\begin{itemize}
\item Scalar invariants: We compute the hyperbolic knot volume and scalar invariants from knot Floer homology including the $\nu$ and $\tau$ invariants, and the Seifert genus.
\item Matrix invariants: We compute the Goeritz matrix of the knot and the bigraded knot Floer homology dimensions (represented as a sparse matrix where entry $(a,m)$ of the Alexander-Maslov grading gives the rank of the homology group with this bigrading).
\item Polynomial invariants: We compute the Alexander and the Jones polynomial and represent them as a zero-padded vector of coefficients. 
\end{itemize}
As a baseline, we also include the braid word, corresponding to the classical student-teacher setup in this case, where the student NN has to learn to reproduce the (PCA-transformed and truncated) embedding learned by the teacher network. For a perfect student, the reconstruction loss from the braid word should be zero.

We train the student NN to reconstruct the embedding for all knot invariants listed above and compare the reconstruction loss of the embedding for different knot invariants as inputs. This provides a measure for whether this invariant is encoded in the embedding data of the teacher. If the embedding learned by the NN employs one of the known invariants, the student should be able to learn the embedding from the given invariant with high accuracy.

The advantage of this setup is that we can also use it to ``subtract off'' the known invariants to uncover new invariants learned by the NN.\footnote{We thank Alex Davies for suggesting this to us.} For each point on the embedding manifold learned by the teacher and presented as a vector in $\mathbbm{R}^n$, we can subtract all the vectors corresponding to maps from the known knot invariants into $\mathbbm{R}^n$ as learned by the student. In the generative case, the linear subtraction is justified by the fact that the embeddings and outputs are related only by a linear output layer. The remaining irreducible parts correspond to information that cannot be constructed from the known knot invariants. We have not carried out a systematic analysis of this but will return to it in the future.

\subsection{Jones Unknot Problem\label{sec:idea_jones_unknot}}
The Jones polynomial is a polynomial knot invariant. It is an open problem whether the Jones Polynomial detects the unknot, i.e., whether 
\begin{align}
J(K)=1 \qquad\Leftrightarrow\qquad K=\text{unknot}\,.
\end{align}
This is the so-called \textit{Jones Unknot Conjecture}, and it is an open problem at the moment. Reasons to believe that the conjecture might be false are that other polynomial invariants (such as the Alexander polynomial) cannot detect the unknot, and that it is known that the Jones polynomial of inequivalent knots can be the same, so why should this not also be the case for the unknot. Reasons to believe that it might be true is that the Jones polynomial can be thought of as a character in bi-graded Khovanov homology~\cite{khovanov2000}, which categorifies the Jones polynomial and is known to detect the unknot~\cite{kronheimer2011khovanov}. Moreover, extensive brute-force searches~\cite{tuzun2018verification} have revealed that there is no other knot up to 24 crossings that also has Jones polynomial $J(K)=1$.

We will use the methods developed here, i.e., generative models from knot invariants to search for counter-examples to this conjecture. In more detail, we will train a NN whose input is the Jones polynomial and whose output is a (randomly chosen) braid word representation of a knot with this Jones polynomial. Given the results of~\cite{tuzun2018verification}, we train the NN with knots that have a large crossing number (and/or complicated knot representations of simple knots that have a large number of crossing).

\begin{figure}[t]
\centering
\subfloat[Jones coefficient embeddings\label{fig:JonesSampling}]{\includegraphics[height=.35\textwidth]{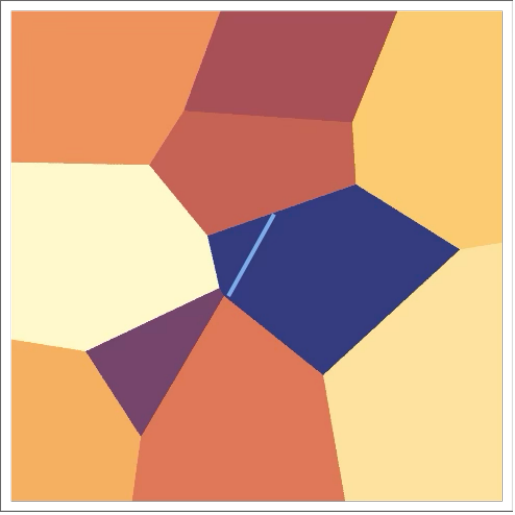}}\qquad\qquad
\subfloat[Sampled paths between knots\label{fig:PathSampling}]{\includegraphics[height=.35\textwidth]{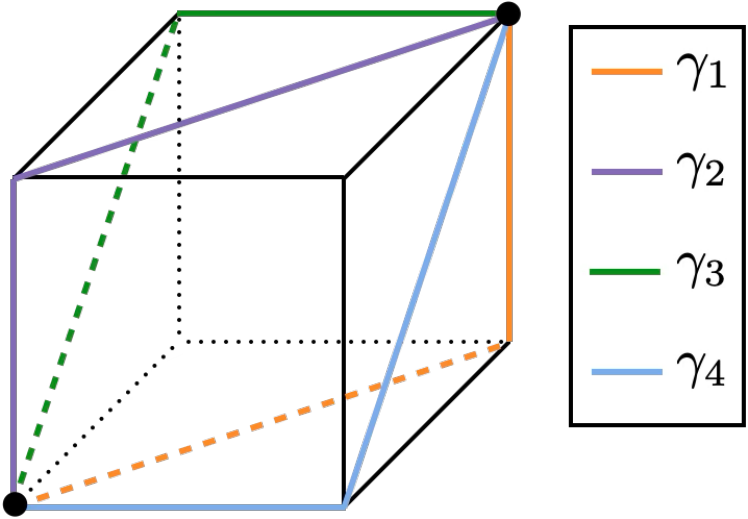}}
\caption{(a) Illustration of sampling in the Jones polynomial embedding space. (b) Different paths sampled between knots, drawn on a cube for better visibility.}
\end{figure}

Equipped with the NN, we can then sample the ``vicinity'' of the Jones polynomial $J(K)=1$ to attempt to find knots that have a trivial Jones polynomial but are not the unknot. While the Jones polynomial coefficients are discrete, they embed into a continuous embedding space. In this space, we associate a random point with the Jones polynomial coefficient with the closest embedded token, where closest is measured with respect to some metric in the embedding space. Choosing the Euclidean metric, we obtain a Voronoi diagram, as usual in unsupervised clustering with decision boundaries that maximize the distance between anchor points (see e.g.~\cite{Ruehle:2020jrk} for an introduction). The idea is illustrated on the left of Figure~\ref{fig:JonesSampling}, where the differently colored regions correspond to the different assignments to Jones polynomial coefficients. Along the boundaries of the region, the Jones polynomial changes.  

Since we use a point in the embedding space as the starting point for our decoder, i.e., our knot generator, the generated braid word changes as we move through the embedding space. However, braid words are also discrete, so there are regions in the embedding space that lead to the same knot. In terms of the partitioning of the embedding manifold according to the knot diagram, we know that if the Jones polynomial changes, so must the knot. Thus, each time we cross over into a differently colored region, signaling that a Jones polynomial coefficient has changed, the decoder will generate a new knot. However, multiple knots can have the same Jones polynomial, and we are looking for such examples where the Jones polynomial is trivial. This means that the colored regions are further subdivided into subregions across which the decoder will generate different (inequivalent) knots, illustrated by the light blue line running through the dark blue region in Figure~\ref{fig:JonesSampling}. We hence want to sample the cell with trivial Jones polynomial in the hope that we discover a substructure along which the knot type changes, thus finding a knot that is not equivalent to the unknot but still has Jones polynomial 1. This would disprove the Jones Unknot Conjecture. 

There are multiple ways to perform this sampling. We explored three possibilities.

\textbf{Gaussian Sampling.} The Jones polynomial is a Laurent polynomial with integer coefficients,
\begin{align}
\label{eq:JonesPolynomial}
J(K) = \sum_i a_i t^i\,,\qquad a_i\in\mathbbm{R}\,.
\end{align}
The Jones polynomial of the unknot, $J(0_1)=1$, corresponds to $a_i=0$ for $i\neq0$ and $a_0=1$ in~\eqref{eq:JonesPolynomial}. Since a generative model will generate a knot for any input (not just integer $a_i$), a simple approach would be to sample $a_i$ ($i\neq0$) from a Gaussian distribution $\mathcal{N}(0,\epsilon)$ for some small variance $\epsilon>0$, and $a_0=1 + \mathcal{N}(0,\epsilon)$. However, since the Transformer operates on tokenized input, this approach is not well-suited. We would have to introduce a new token for each random, real coefficient, leading to completely out-of-distribution data. Instead, we can add Gaussian noise to the token embeddings, which is what we opt for here. More precisely, we add a Gaussian noise sampled i.i.d.\ from $\mathcal{N}(0,\epsilon)$ to the output of the encoder network. As we increase $\epsilon$, the NN will eventually generate knots that are not in the equivalence class of the unknot, and the hope is that this happens for small enough $\epsilon>0$ such that the Jones polynomial of this knot is still 0. 

\textbf{Trajectory Sampling.} The idea is similar to the Gaussian sampling procedure described above. However, instead of injecting random Gaussian noise after the encoder, we sample trajectories in the embedding space from the embedded Jones polynomial of the unknot to the embedded Jones polynomial of other knots with simple Jones polynomials.\footnote{We define a Jones polynomial to be simple if its span, which is defined as the difference between the largest power and the smallest power that occurs in the polynomial, is small.} If knots with simple Jones Polynomial cluster in embedding space, the hope is we can find more knots along this trajectory with simple Jones polynomials, and that among them is a knot with Jones Polynomial 1. Let us denote the vectors corresponding to the embedded tokens of the zero-padded Jones polynomial of the unknot and the other simple knot $\vec{v}_{i}$ and $\vec{w}_{i}$, respectively. Here, $i$ runs over the input length, which is the number of components of the zero-padded vector of Jones Polynomial coefficients $a_i$ in~\eqref{eq:JonesPolynomial}. We construct four different types of piecewise linear trajectories along which we sample 10 equidistant points:
\begin{itemize}
\item We interpolate linearly the $a^\text{th}$ coordinate of all vectors $\vec{v}_{i}$ to the $a^\text{th}$ coordinate of all vectors $\vec{w}_{i}$. Subsequently, we linearly interpolate all remaining coordinates of $\vec{v}_{i}$ to those of $\vec{w}_{i}$. Schematically we write such paths as $\gamma_1(a)$ and e.g.\ $\gamma_1(1) = (t,0,\ldots,0)\circ(0,t,\ldots,t)$ for the first component, where $\circ$ means we do the LHS action first, then the RHS action.
\item Same as above, but we invert the order. We first linearly interpolate all but the $a^\text{th}$ coordinate from $\vec{v}_{i}$ to $\vec{w}_{i}$, followed by interpolating the missing $a^\text{th}$ component. Schematically, $\gamma_2(1) = (0,t,\ldots,t)\circ(t,0,\ldots,0)$ for $a=1$.
\item Similar to the first case, but we simultaneously linearly interpolate $2$ coordinates at the $a^\text{th}$ and $b^\text{th}$ position, followed by linearly interpolating the remaining $N-2$ components, $\gamma_3(1,2) = (t,t,0,\ldots,0)\circ(0,0,t,\ldots,t)$ for $a=1$, $b=2$.
\item Same as above, but we invert the order, so we get schematically the paths $\gamma_4(1,2) = (0,0,t,\ldots,t)\circ(t,t,0,\ldots,0)$ for $a=1$, $b=2$.
\end{itemize}
We illustrate the different trajectories in $\mathbbm{R}^3$ in Figure~\ref{fig:PathSampling}.

\textbf{Beam and Temperature Sampling.} Instead of injecting noise at the encoder level, we can also inject noise at the embedding layer, i.e., the last layer before output of the decoder. There are multiple ways to do this, and we chose to combine temperature sampling together with a beam search. Adding temperature means that instead of greedily selecting the most probable letter as predicted by the final softmax layer of the decoder network, we sample the logits $\ell_i$ from a probability distribution with a non-zero temperature, yielding logit probabilities 
\begin{align}
p_i = \text{softmax}(\ell_i/T).
\end{align}
Note that in the $T\to0$ limit, this selects the largest logit and thus is equivalent to greedy sampling. Instead of sampling all braid words independently this way, we combine this with a beam search, where we sample multiple times and keep the $N$ braid words with the highest combined score. We choose a beam of size 100 for this experiment.

\section{Experimental Results\label{sec:expts}}

In this section we present the results of experiments probing the different learning setups described in Section~\ref{sec:ideas}. 

\subsection{Knot Data\label{sec:exp_data}}
In this section we describe the algorithms used to generate representatives of knot classes. We begin with the algorithm \textsc{RandomKnot}, Algorithm 6 of \cite{gukov2020learningunknot}. This algorithm takes $n_\text{letters}$, $n_\text{strands}$, and $M$ as arguments. It first generates a random braid with these parameters from the uniform distribution on the Artin generators, each represented by $\pm k$, where \nobreak{$k = 1,\dots,n_\text{strands}-1$}. The random braid in general gives a link, but may be turned into a knot via the \textsc{Knotify} algorithm of \cite{gukov2020learningunknot}. The algorithm then applies $M$ random Markov moves and braid relations to the braid, and uses another algorithm that makes the braid smaller while preserving its class. We take $M=1$, since further scrambles will occur later. This procedure may produce a braid of length different from $n_\text{letters}$, in which case the procedure is repeated until the braid has the desired length. The output is a random knot, represented as a braid word.

Next we take this seed braid word and generate different representatives of the same class. We do so by applying $n_\text{scrambles}$ random Markov moves to the seed braid word.  This yields a list of braid words, potentially of different length due to the scrambling, each a different representative of the same knot class. For training of the auto-regressive generator, we use the shortest braid word from each class as targets. After identifying this representative, the braids are all padded to the same length via stabilization Markov moves. In particular this means that datasets with a fixed $n_\text{letters}$ may have variance in $n_\text{strands}$ due to the stabilizations.

Finally, we perform an ablation study, for which we generate data across a wide range of parameters. For each training dataset we generate $50,000$ equivalence classes of braids with $20$ representatives each, with $5$ strand braid words (prior to stabilization). We keep the number of strands fixed but ablate
\begin{equation}
\label{eq:train_data_param_variance}
n_\text{letters} \in \{20,25,30,35,40,45,50\}\,, \qquad n_\text{scrambles} \in \{5,10,15,20\}\,,
\end{equation}
since we intuitively expect that the braid length and number of Markov moves performed will have the greatest impact on learning. For test / validation datasets, we use the same parameters but take the number of equivalence classes to be $25$. This yields about $10$GB of data, with train datasets of size about $150$MB to $400$MB as $n_\text{letters}$ varies from $20$ to $50$.

\subsection{Contrastive Invariance Learning\label{sec:exp_contrastive_invariance}}
\subsubsection{Experiments}
As explained in Section~\ref{sec:idea_invariance}, we tried standard feed-forward NNs and CNNs. The latter architecture is constructed to be invariant under braid conjugations. We performed some mild hyperparameter tuning for both the semi-hard triplet loss and the centroid loss. The number of training epochs has the biggest impact, and the models tend to overfit. As preventative measures we tried weight regularization and dropout layers, but the best results were obtained by early stopping. Also, the losses work better with larger batch sizes (we used 128), and learning rates can be rather large. As an example, we will illustrate results here using the dataset of 50k knots with braid length 30 with 5 strands, 20 knots per equivalence class, and 10 scrambling moves. Very small networks give very good results, and we used the following:
\begin{itemize}
\item \textbf{MLP:} We choose 2 hidden layers with 64 nodes and tanh activation (7k parameters), 300 epochs, and Adam optimizer with learning rate $10^{-3}$.
\item \textbf{CNN:} We choose a 1D convolutional layer with 64 filters and leaky ReLU, followed by a two-layer MLP with 64 nodes and tanh activation (130k parameters), 500 epochs, and Adam optimizer with learning rate $10^{-4}$.
\end{itemize}
Despite using a 16-dimensional embedding layer, a PCA reveals that the trained NNs only use 1 or 2 embedding dimensions. Because of this, we do not even perform a PCA and plot the first two embedding dimensions directly. For better visibility, we apply a standard scaler, i.e., we map the outputs to zero mean and unit variance in the embedding space.

\begin{figure}[tb]
\centering
\includegraphics[width=.43\textwidth]{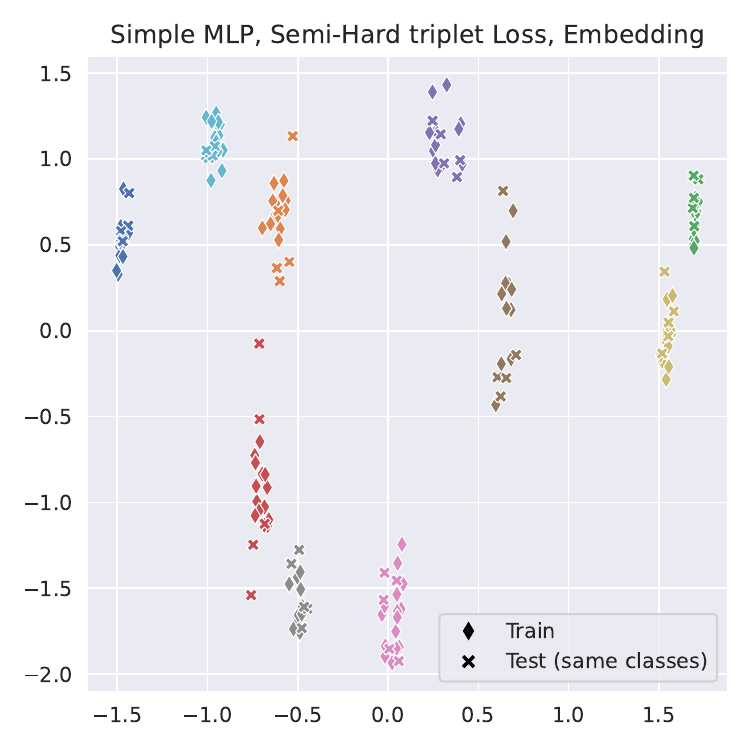}\qquad
\includegraphics[width=.43\textwidth]{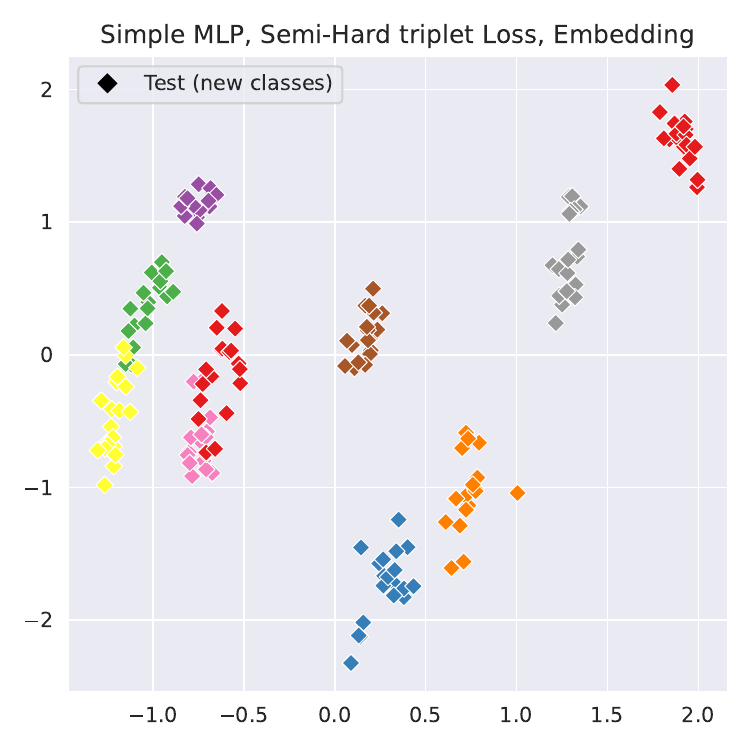}\\
\includegraphics[width=.43\textwidth]{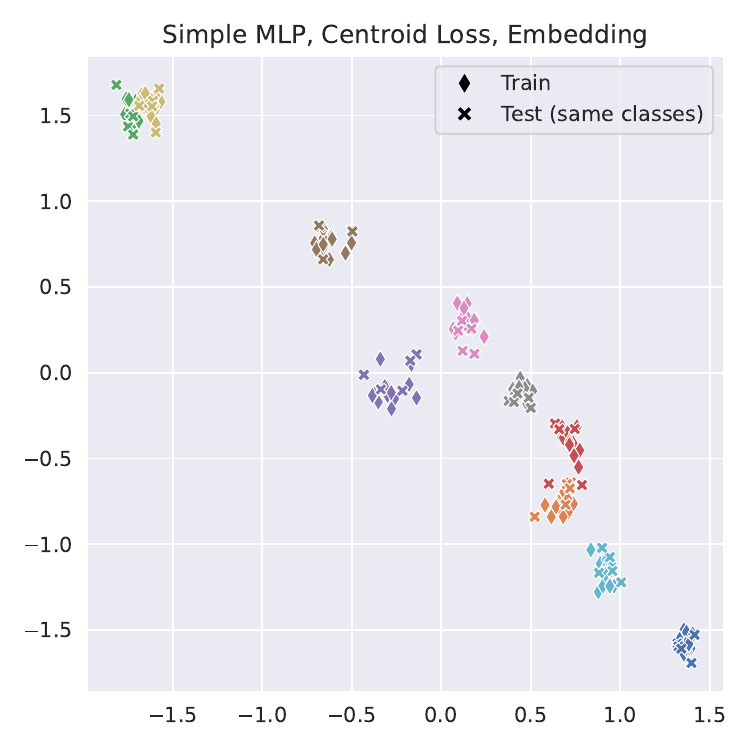}\qquad
\includegraphics[width=.43\textwidth]{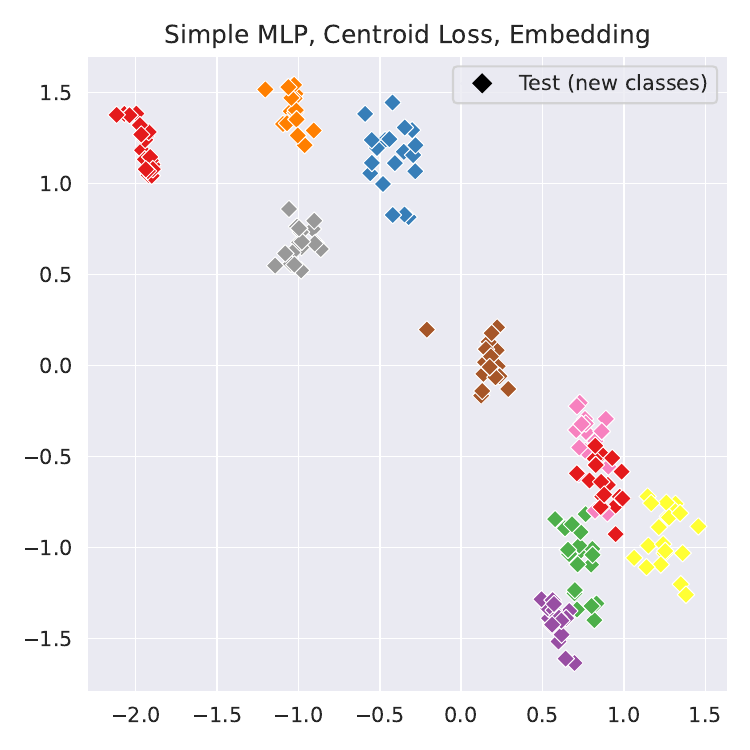}\\
\caption{Semi-hard triplet loss and centroid loss for the MLP Networks.}
\label{fig:ContrastiveLearningClustering1}
\end{figure}

\begin{figure}[tb]
\centering
\includegraphics[width=.43\textwidth]{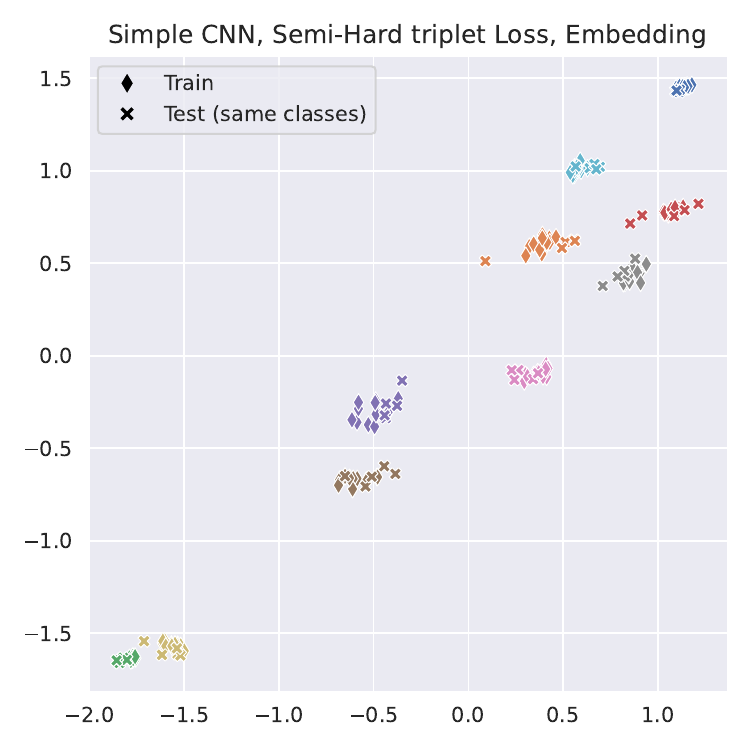}\qquad
\includegraphics[width=.43\textwidth]{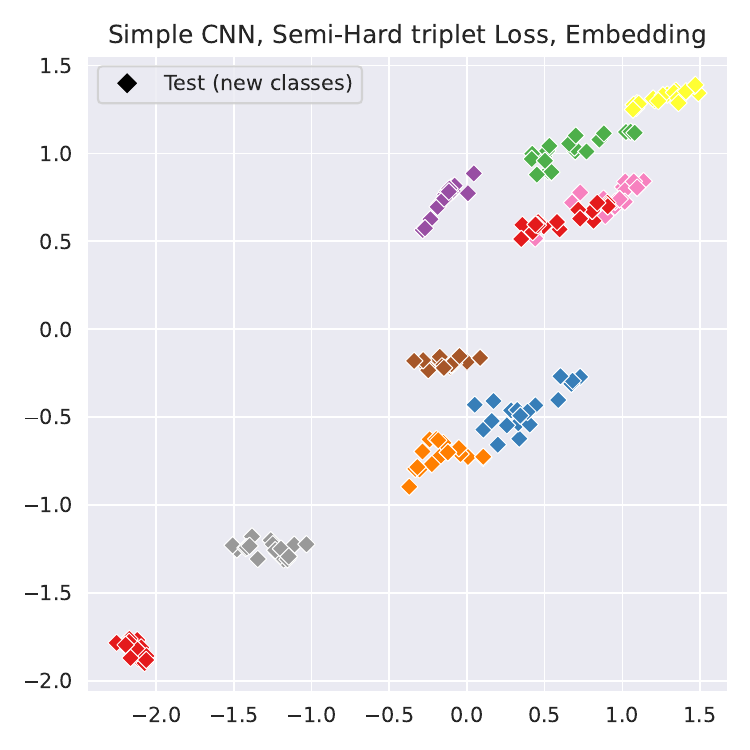}\\
\includegraphics[width=.43\textwidth]{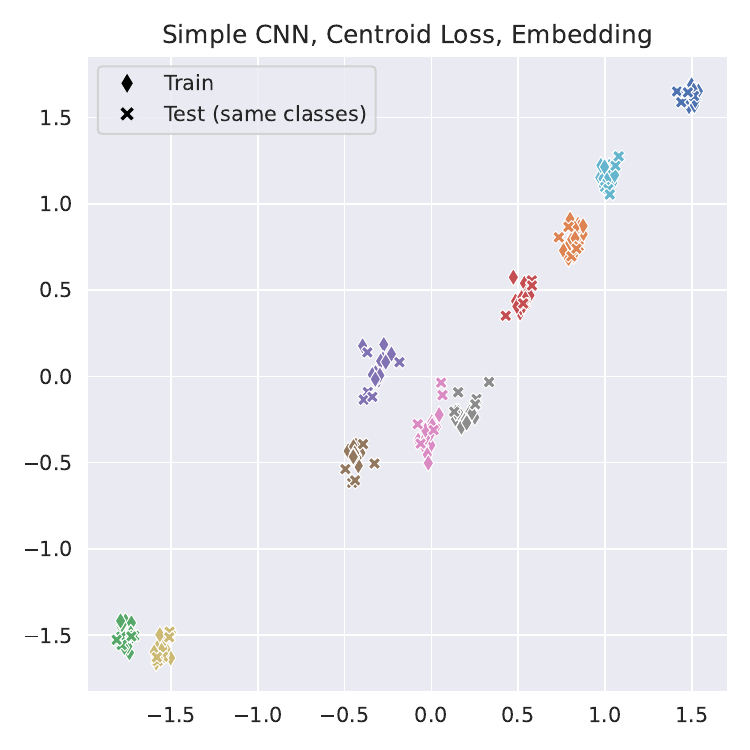}\qquad
\includegraphics[width=.43\textwidth]{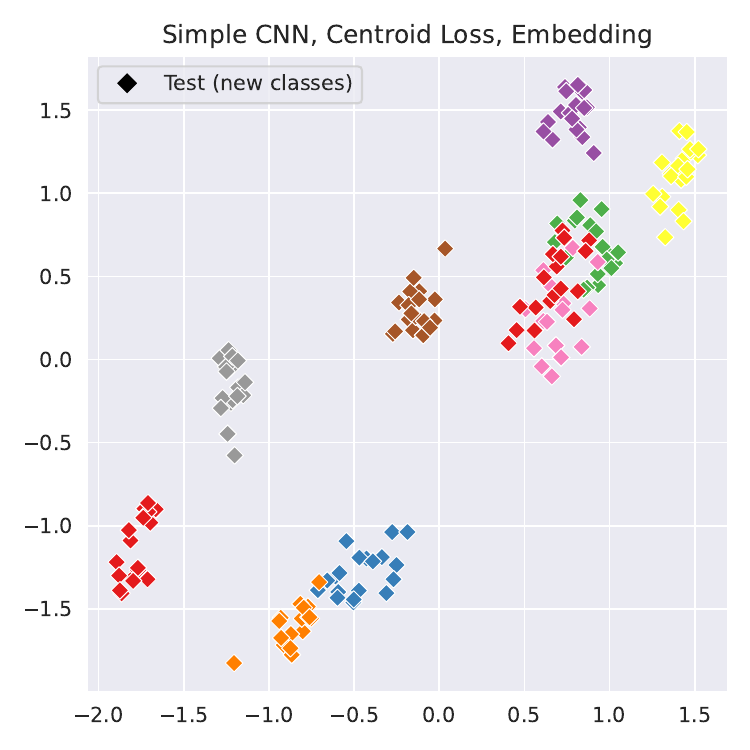}
\caption{Semi-hard triplet loss and centroid loss for the CNN Networks.}
\label{fig:ContrastiveLearningClustering2}
\end{figure}

There are two conceptually inequivalent ways of generating a test set.
\begin{itemize}
\item \textbf{In-distribution:} We take inequivalent knot representatives which, however, belong to a knot equivalence class that the NN has been trained with
\item \textbf{Out-of-distribution:} We use knot representatives from equivalence classes which the NN has not seen during training.
\end{itemize}
To generate these, we only train the NNs with 15 knot representatives per equivalence class, keeping 5 for the `in-distribution' test set. Similarly, we do not train it with all 50k equivalence classes but keep 50 percent for the `out-of-distribution' test set. We find that in-distribution generalization works very well, meaning that the test set knots get mapped to the corresponding (tight) cluster of training knots. For the out-of-distribution knots, the clustering is a bit less tight. To illustrate this, we plot the results for 10 in-distribution and 10 out-of-distribution knot equivalence classes for the MLP and CNN Figures~\ref{fig:ContrastiveLearningClustering1} and~\ref{fig:ContrastiveLearningClustering2}, respectively. For the in-distribution test set, we also include the 15 equivalent knots per cluster from the training set. For the out-of-distribution test knots, we plot all 20 equivalent knots.

\subsection{Generative Invariance Learning\label{sec:exp_generative_invariance}}

Next, we use the generative method for learning topological invariance of knots, as described in Section~\ref{sec:idea_invariance}. We begin with an ablation study (including transfer learning and beam search) that sets the stage for further analyses, and then proceed to PCA and clustering analyses in the embedding dimension, and longer training runs. 

\subsubsection{Ablation Study}
We train our models with the Transformer library Int2Int~\cite{f_charton_Int2Int}, for all train data sets with different values of $n_\text{letters}$ and $n_\text{scrambles}$ given in \eqref{eq:train_data_param_variance}. We train with batch size $64$ for $25$ epochs, where each epoch trains on $300,000$ knots. 
The encoder and decoder embedding layers each have 256 dimensions. Both the encoder and decoder consist of positional embedding and 4 Transformer layers, with  8 attention heads per layer. Each encoder and decoder layer includes 1 feed-forward network layer by default, with \texttt{ReLU} activation.

During training, validation is performed every $200$ batches, yielding about $23$ validations per epoch. Each time, the decoder is used to generate a new knot $K'$ using the embedding of a test set knot $K$. If the decoder is well-trained, $K'$ and $K$ should be equivalent knots. To perform validation on $K$ and $K'$, one wishes to utilize a fast-to-compute topological invariant that is sufficiently robust (at the large crossing numbers we are considering) to distinguish between equivalent and non-equivalent knots. For this purpose, we utilize the hyperbolic volume $\text{vol}(K)$ of the knot complement, computed in snappy, where we consider two volumes to be equivalent if $|\text{vol}(K)-\text{vol}(K')|\leq 10^{-5}$.

\begin{figure}[t]
\centering
\begin{tabular}{cc}
\includegraphics[width=0.4\textwidth]{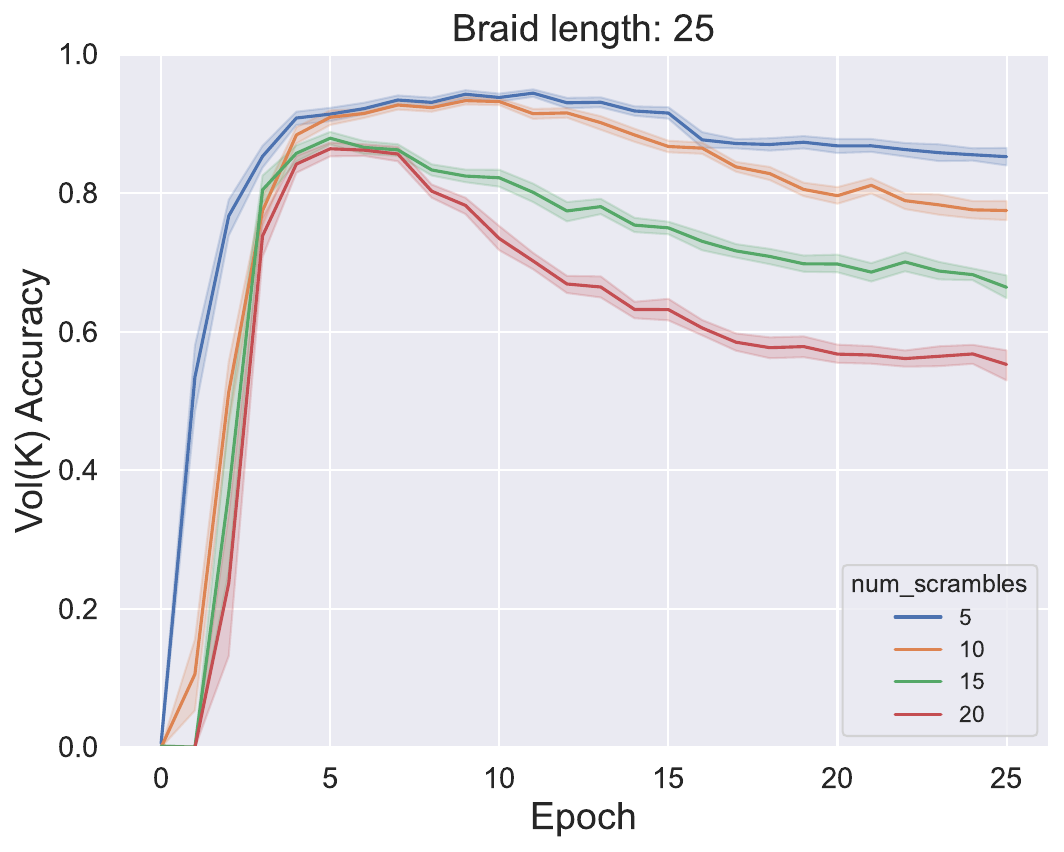} &
\includegraphics[width=0.4\textwidth]{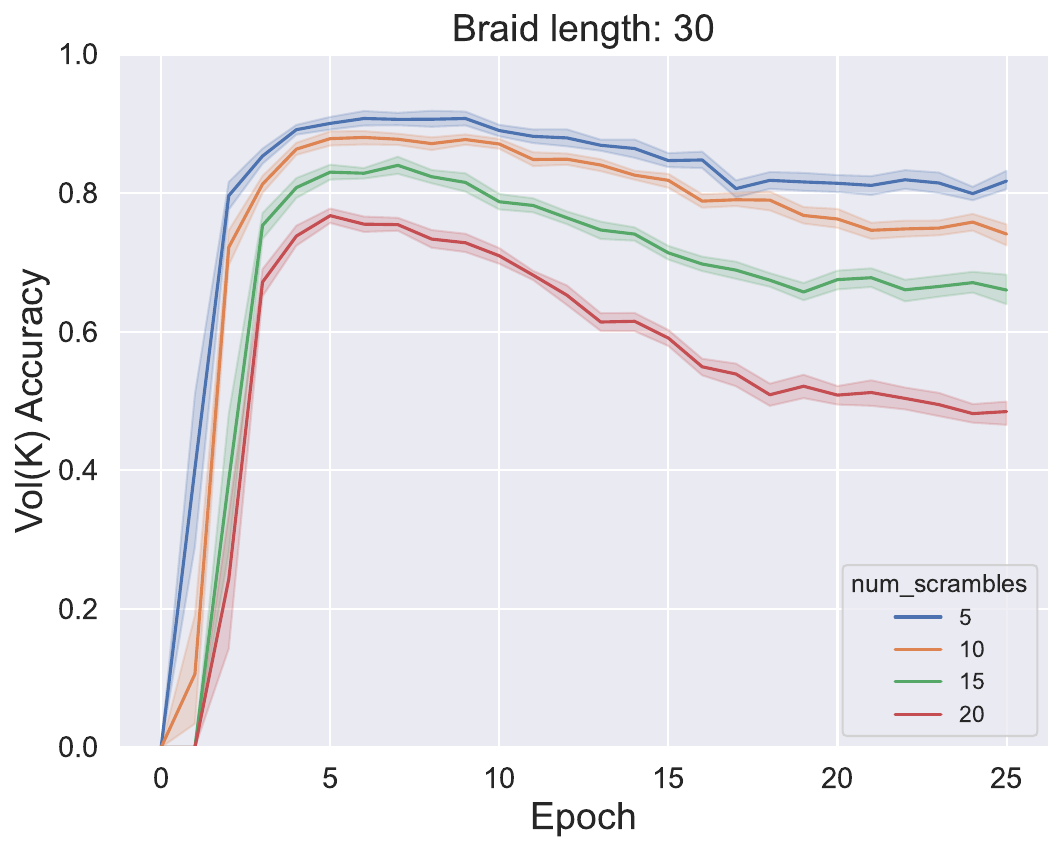} \\
\includegraphics[width=0.4\textwidth]{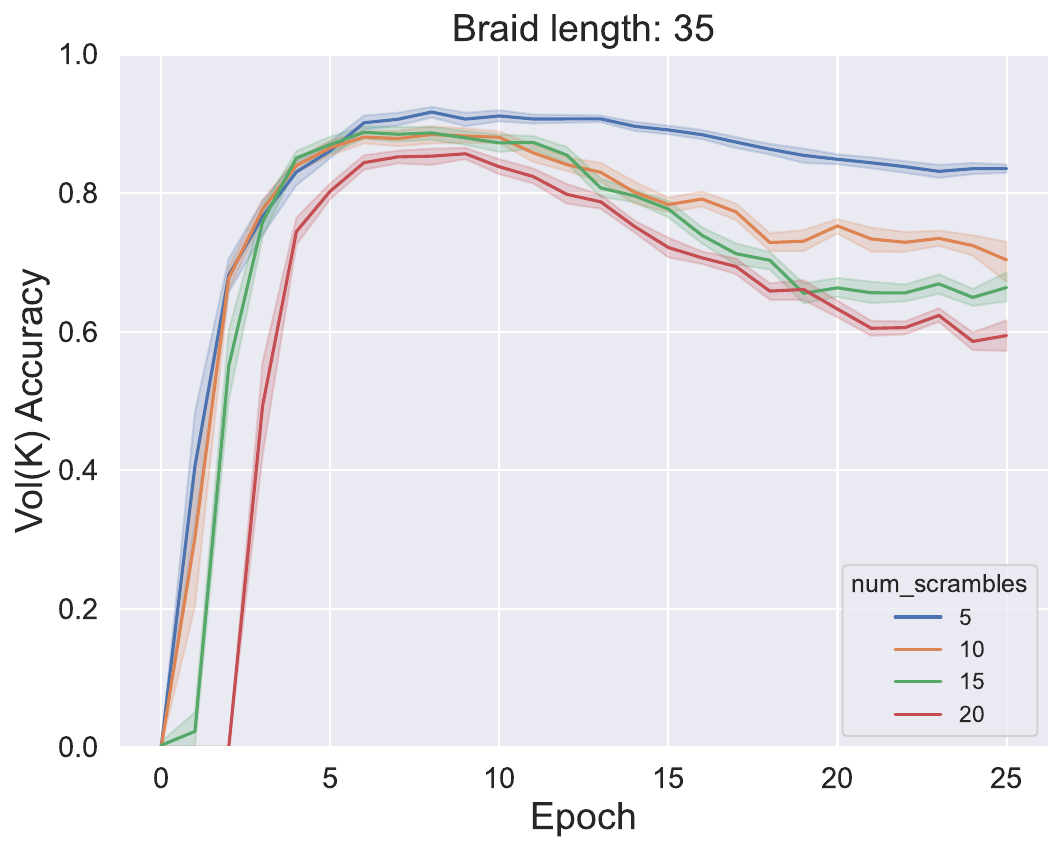} &
\includegraphics[width=0.4\textwidth]{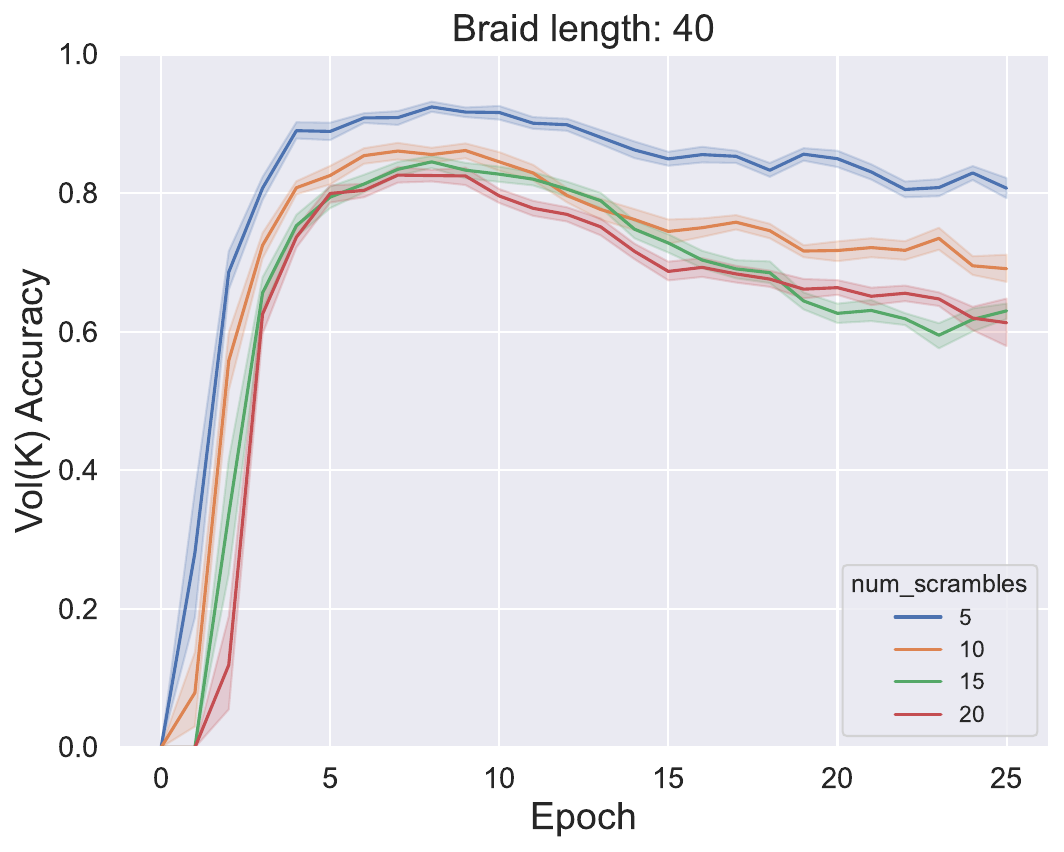} \\
\end{tabular}
\caption{Ablation studies for different braid lengths. The $\text{vol}(K)$ test accuracy decreases monotonically with increasing number of scrambles.}
\label{fig:ablation_braid_lengths_fixed}
\end{figure}

\begin{figure}[ht]
\centering
\begin{tabular}{cc}
\includegraphics[width=0.4\textwidth]{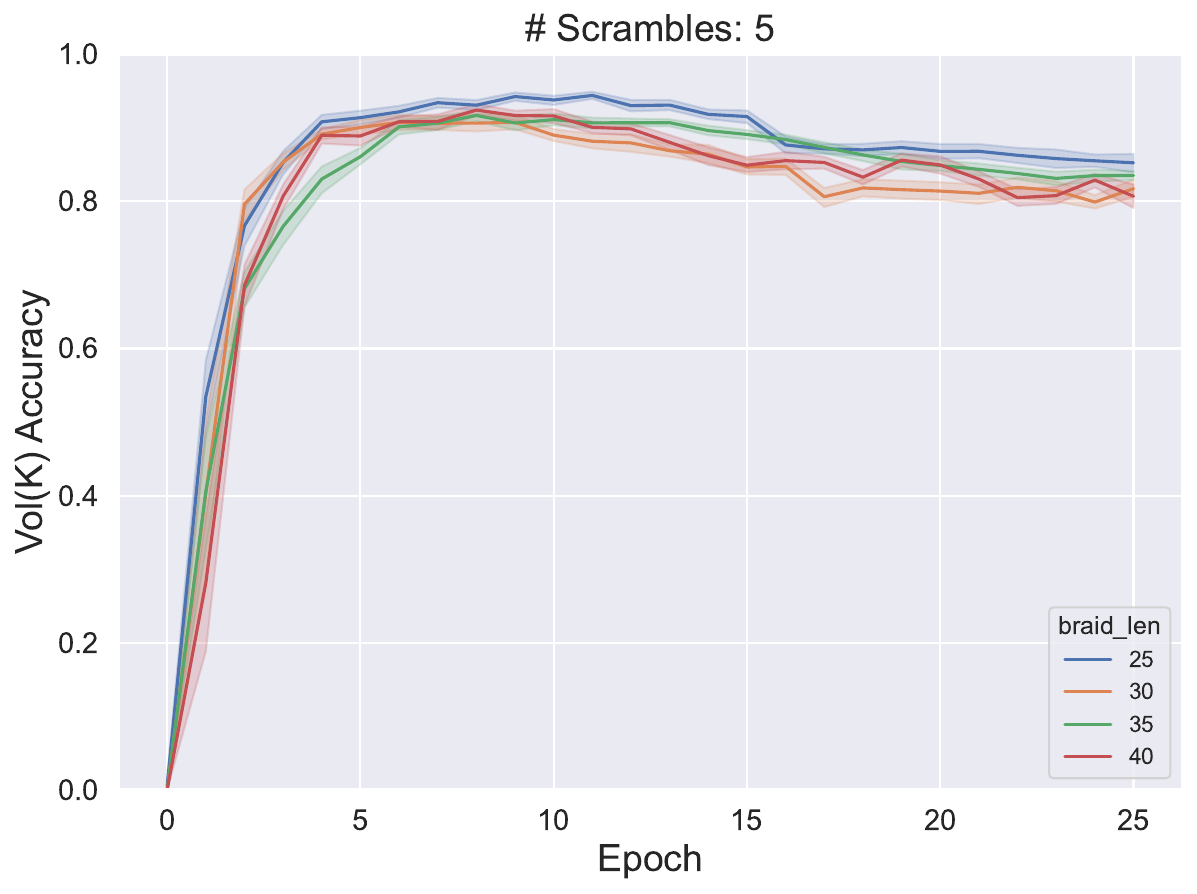} &
\includegraphics[width=0.4\textwidth]{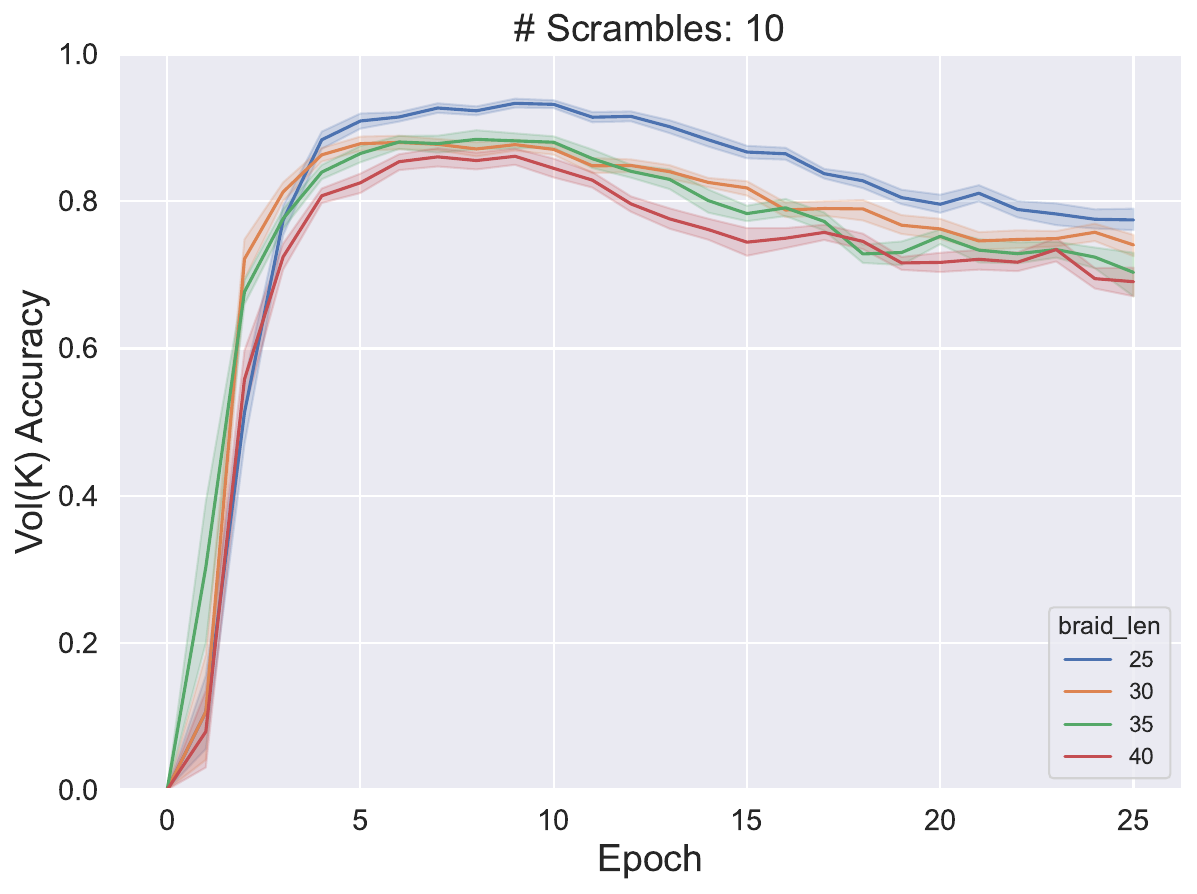} \\
\includegraphics[width=0.4\textwidth]{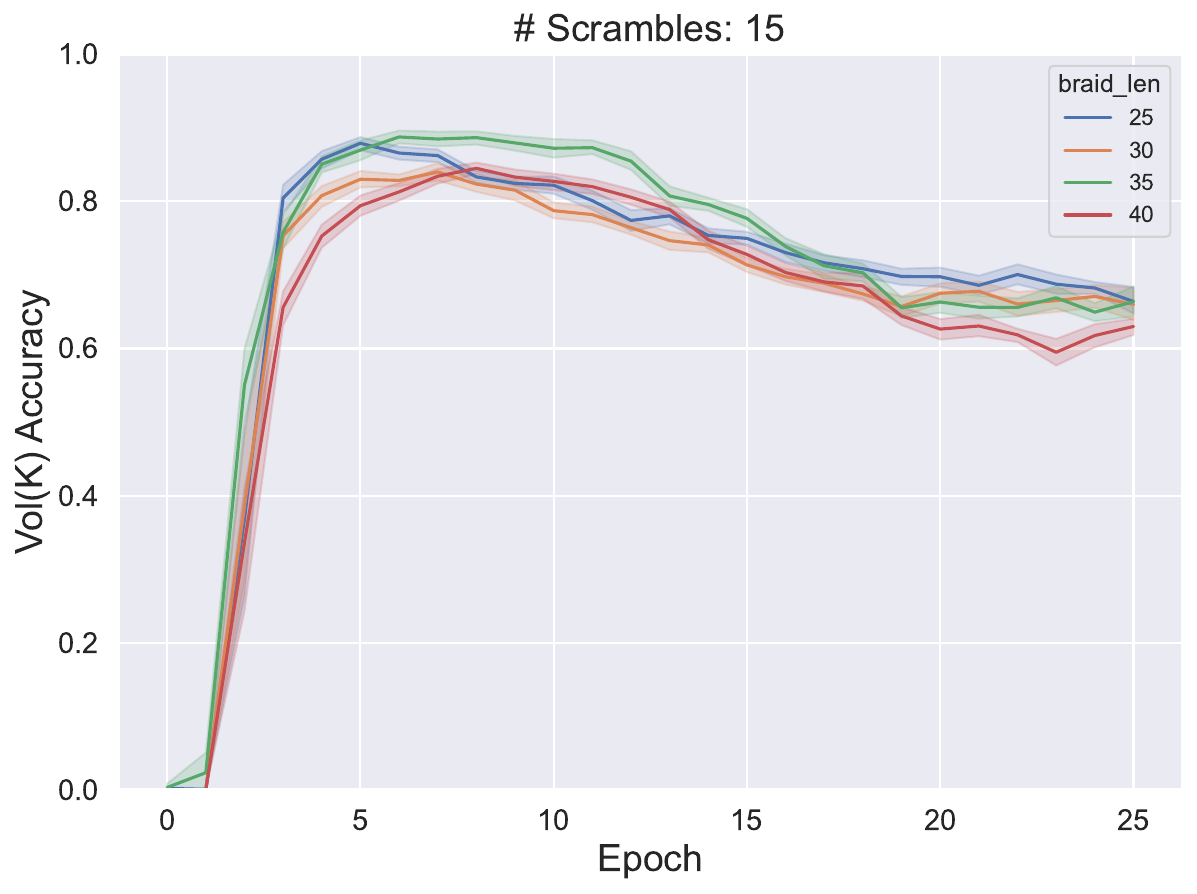} &
\includegraphics[width=0.4\textwidth]{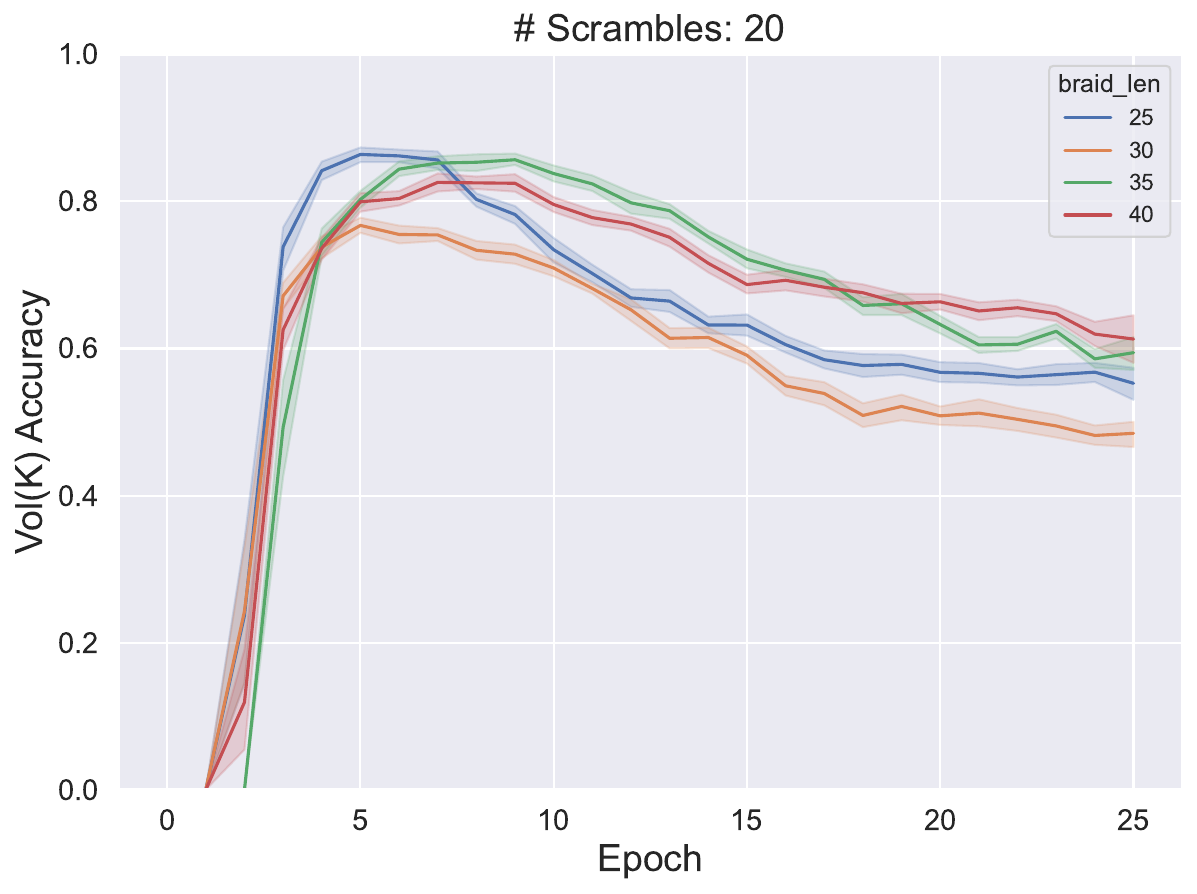} \\
\end{tabular}
\caption{Ablation studies for different numbers of scrambles. The $\text{vol}(K)$ test accuracy does not demonstrate a clear correlation with the braid length.}
\label{fig:ablation_num_scrambles_fixed}
\end{figure}

The results of the study are presented in Figures~\ref{fig:ablation_braid_lengths_fixed} and~\ref{fig:ablation_num_scrambles_fixed}. In Figure \ref{fig:ablation_braid_lengths_fixed}, each plot has a fixed braid length and a varying number of scrambles. We see that the $\text{vol}(K)$ test accuracy decreases monotonically with increasing number of scrambles. In Table \ref{fig:ablation_num_scrambles_fixed}, each plot has a fixed number of scrambles and a varying braid length. We see that the $\text{vol}(K)$ test accuracy does not demonstrate a clear correlation with the braid length, showing that the Transformer performs well across the entire range of braid lengths considered.

We have also performed a PCA and clustering analysis of the embeddings of the train and test knots, for $n_\text{letters}=30$ and $n_\text{scrambles}=10$. Interestingly, in this case the embeddings of knots produced by a randomly initialized Transformer were well-clustered by knot class, with the first two principal components making up about $50\%$ of the explained variance. In the limit that $n_\text{scrambles}$ go to zero, this is a natural expectation, and also that the untrained clustering would persist for small enough $n_\text{scrambles}$. This explains the clustering result. The untrained network (of course) had zero accuracy on $\text{vol}(K)$ prediction: the network embeds the braids in close proximity to one another due to similar braid words, but knows nothing of their topology. In contrast, the trained Transformer spreads the knots out in the embedding dimension, with the first two principal components making up about $11\%$ of the explained variance. The $\text{vol}(K)$ test accuracy was $92.4\%$ for the trained network.

In Table \ref{tab:beam_search_results} in Appendix~\ref{app:BeamSearch}, we present the results of a beam search with the Transformer, trained on knots with $n_\text{letters}=30$ and $n_\text{scrambles} = 10$. A fixed knot is input into the encoder, and its embeddings are used by the decoder in the beam search, with the goal of using the generative model to produce many topologically equivalent knots. We see from the Table that all knots from the beam search have the same $\text{vol}(K)$, and therefore we expect them to be topologically equivalent. This equivalence is demonstrated explicitly, where it is shown how various braids may be obtained from one another by conjugation and braid relations.

Finally, we also studied transfer learning, specifically whether a Transformer trained with $n_\text{scrambles}=5$ could generate equivalent knots when the test knot has been scrambled more that the training data. Table~\ref{tab:transfer_learning} presents the associated results, which demonstrate that the trained networks are able to generalize to knots with more scrambles. The  accuracy is $5\%-10\%$ lower than for the test data at $n_\text{scrambles}=5$, demonstrating a reasonable degree of transfer learning.

\begin{table}[t]
    \centering
\begin{tabular}{|c|c|c|c|c|c|c|}
    \hline
   $n_\text{letters}$ & 25 & 30 & 35 & 40 & 45 & 50 \\
\hline \hline
 Train Vol(K) accuracy (\%) & 97.8 & 98.8 & 96.2 & 97.6 & 97.0 & 94.4 \\
 Transfer Vol(K) accuracy (\%) & 95.7 & 89.6 & 96.3 & 91.4 & 90.3 & 89.0 \\ \hline
\end{tabular}
\caption{Transfer learning results. The first row shows the $\text{vol}(K)$ accuracy of the model trained with $n_\text{scrambles}=5$, evaluated on its test data. The second row shows the average $\text{vol}(K)$ accuracy of the same trained model on test data with $n_\text{scrambles} \in \{10,15,20\}$.}
\label{tab:transfer_learning}
\end{table}

\subsubsection{Long Runs}
In our ablation study, in Figures \ref{fig:ablation_braid_lengths_fixed} and \ref{fig:ablation_num_scrambles_fixed}, we see that $\text{vol}(K)$ accuracy decreases after some amount of time. For this reason, we performed a longer training run to see if this was corrected at later times; if not, we would like to understand the overtraining.

\begin{figure}
    \centering
    \includegraphics[width=0.45\textwidth]{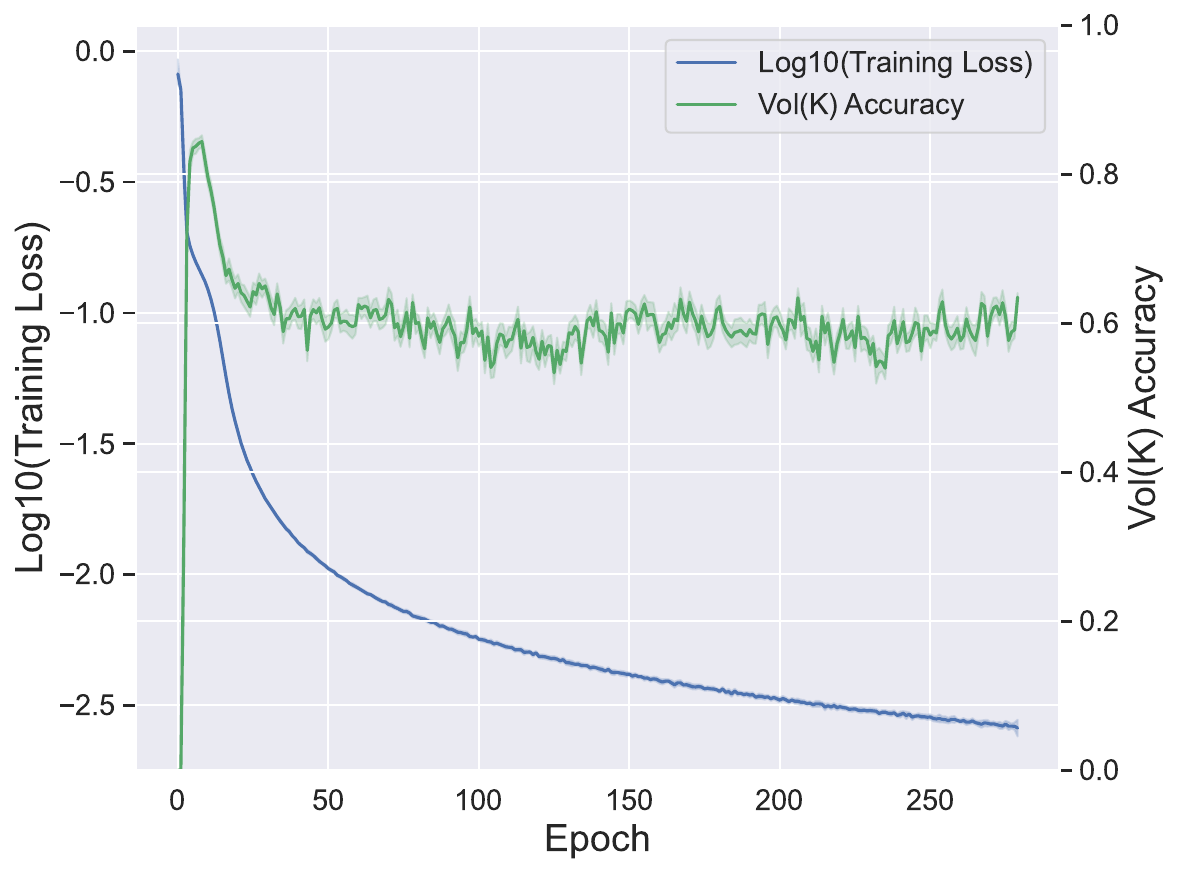} 
    \includegraphics[width=0.45\textwidth]{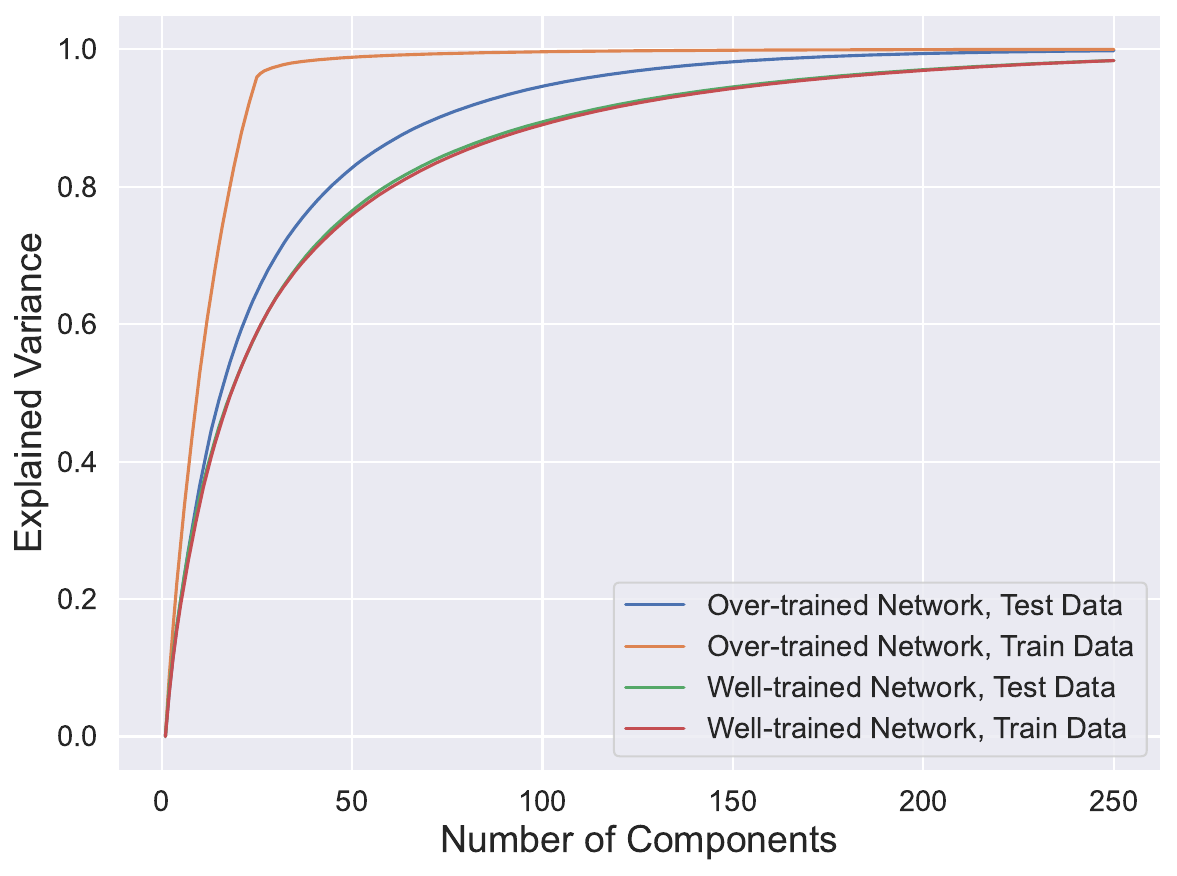} \\ \vspace{.25cm}
    \includegraphics[width=0.9\textwidth]{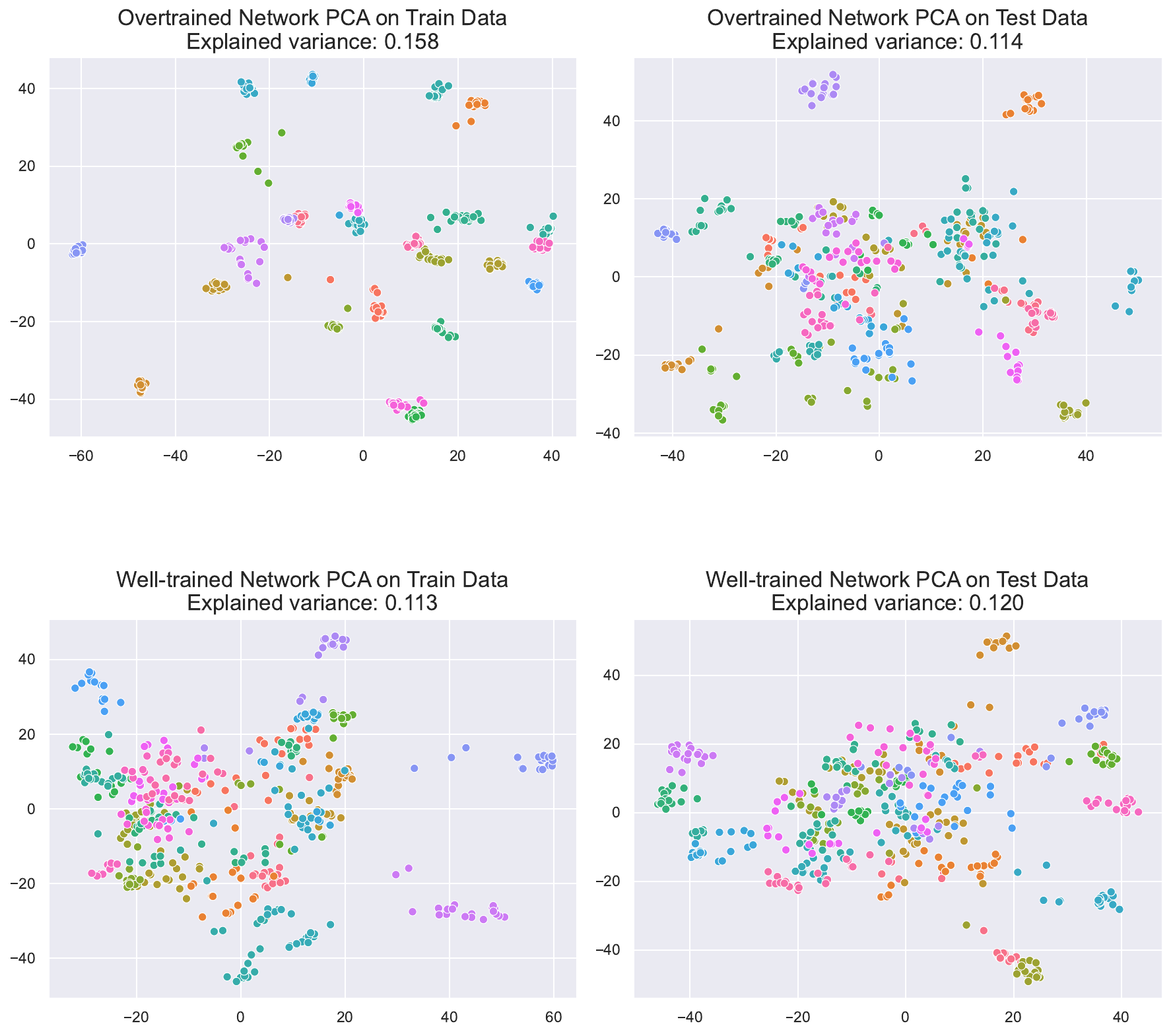}
    \caption{\emph{(Top LHS):} The train loss decreases steadily even though $\text{vol}(K)$ accuracy is maximized early. The network with minimum train loss is over-trained, and the network with maximal $\text{vol}(K)$ accuracy is well-trained. \emph{(Top RHS):} the over-trained network exhibits much higher explained variance with fewer principal components. \emph{(Bottom Two Rows):} the over-trained network is more clustered in the first two PCs than the well-trained network. }
    \label{fig:long_run}
\end{figure}

Specifically, we trained a model with $n_\text{letters}=30$ and $n_\text{scrambles}=15$ for over $250$ epochs. The model architecture and remaining hyperparameters are the same as in the ablation study. The results are shown in Figure \ref{fig:long_run}. We study the behavior of the network at two particular time slices during training. One is at the maximum of $\text{vol}(K)$ accuracy, which we call the well-trained network. The other is at the minimum of the (binary cross-entropy) train loss, which we call the over-trained network. We see from the top LHS plot that the loss flattens slightly as the $\text{vol}(K)$ accuracy is peaking, but then dips strongly as the accuracy leaves its maximum. During most of training, $\text{vol}(K)$ accuracy does not improve even though the train loss decreases steadily. On the top RHS, and in the bottom four plots, we study a principal component analysis of the test and train embeddings for both the over-trained and well-trained network. In the cluster plots of the first two principal components, colors indicate topologically equivalent knots. The over-trained network exhibits tighter clustering than the well-trained network in those plots for both the test and train data. This observation is born out in the explained variance plot, where we see that the well-trained network has nearly identical explained variance on both the test and train data, both of which exhibit significantly less explained variance than for the over-trained network. We see that the highest explained variance occurs for the train data of the over-trained network, indicating that the network has overfit to the train data.

We see that training for a longer time leads to overtraining that hurts performance, hence justifying our choice in the ablation study to study the network parameters at relatively early times where the $\text{vol}(K)$ accuracy is maximized. A similar overtraining was observed in the contrastive learning setup.

\subsection{Extracting Invariants\label{sec:exp_invariants}}
In order to see which known knot invariants are learned by the NNs, and to extract the unknown ones, we use a teacher-student setup as described in Section~\ref{sec:idea_invariants}. 

\subsubsection{Extracting invariants from contrastive learning}
\begin{figure}[t]
\centering
\includegraphics[width=.92\textwidth]{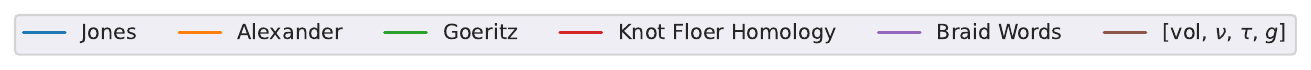}\\[-.15cm]
\includegraphics[width=.44\textwidth]{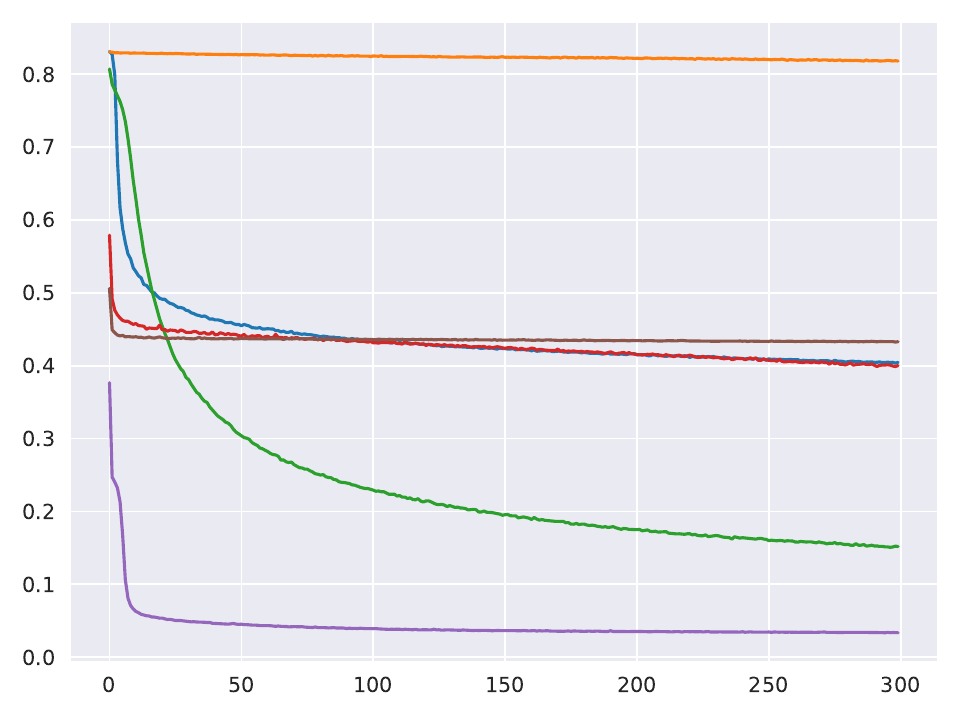}\qquad
\includegraphics[width=.44\textwidth]{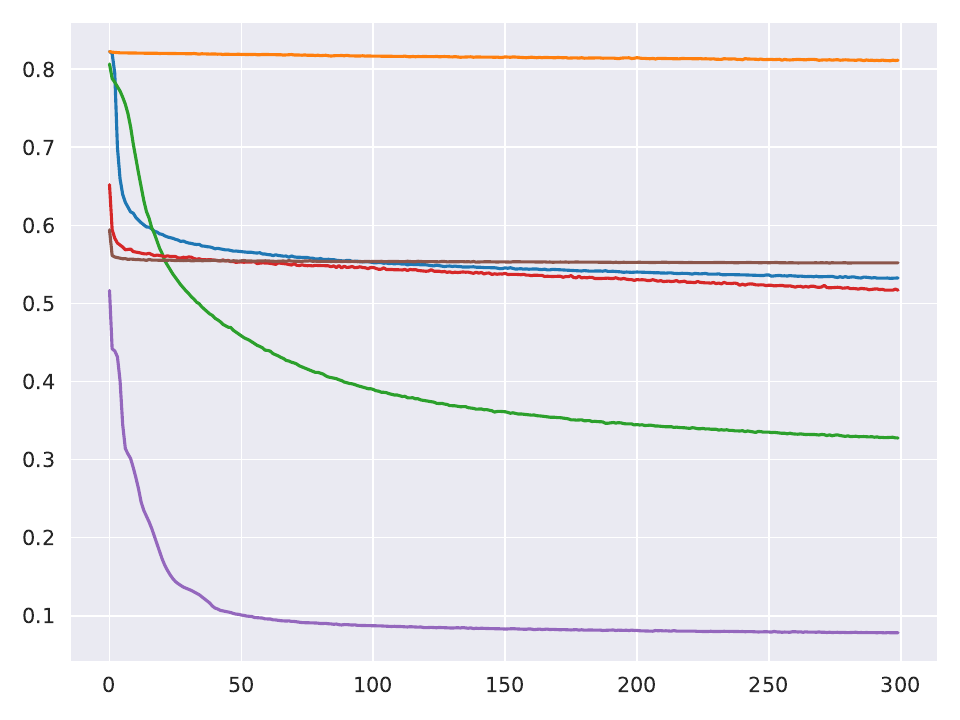}\\
\includegraphics[width=.44\textwidth]{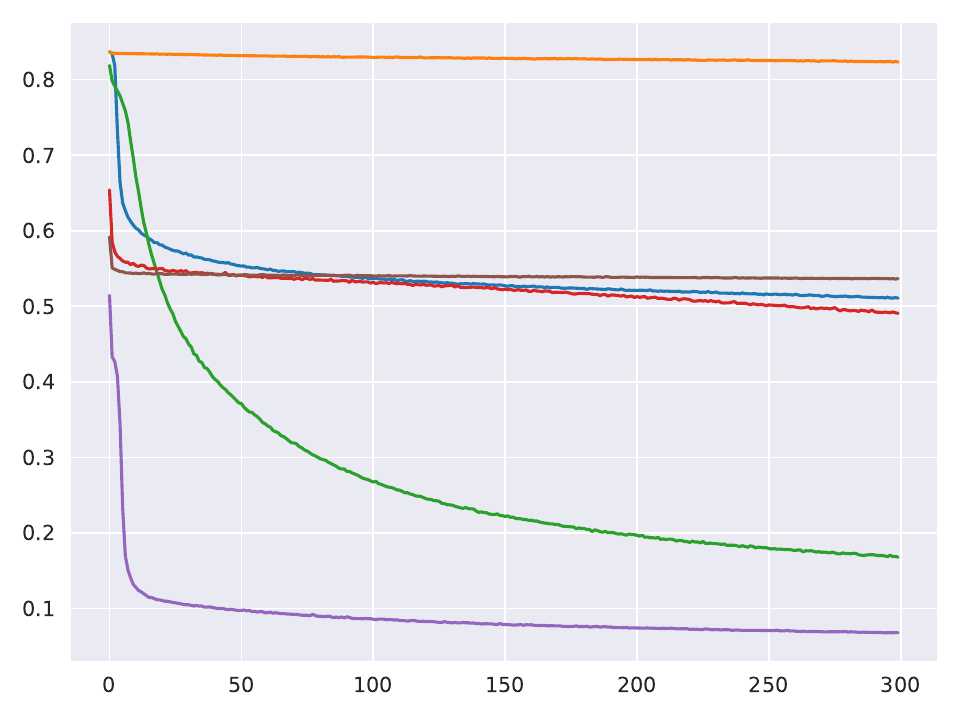}\qquad
\includegraphics[width=.44\textwidth]{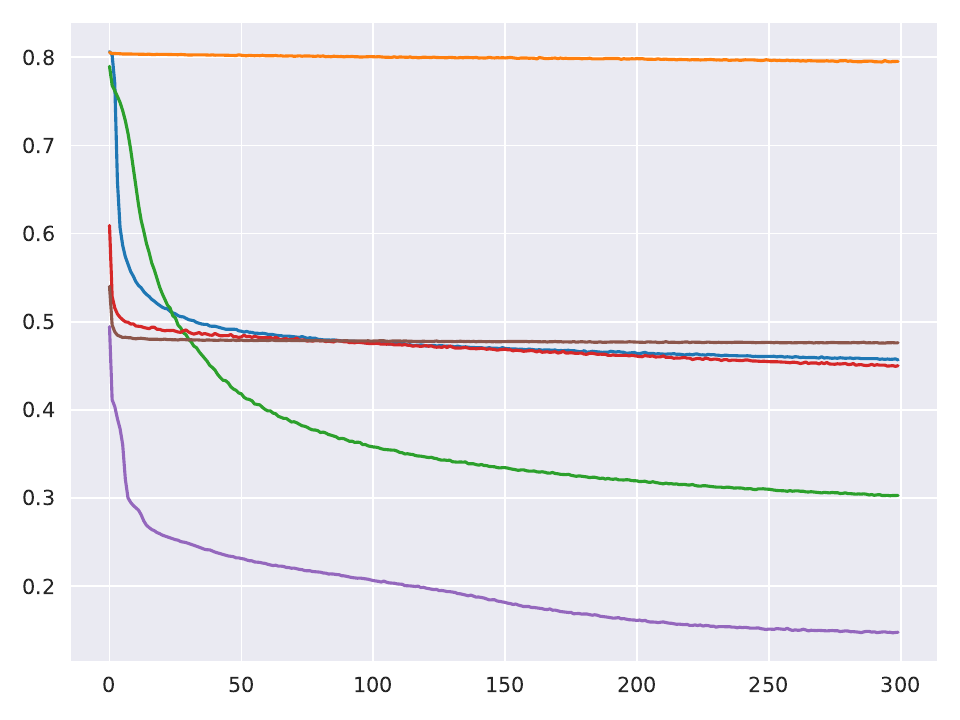}
\caption{Invariants learned by the Teacher NN trained with contrastive loss. Top Left: MLP with semi-hard triplet loss. Top Right: MLP with centroid loss. Bottom Left: CNN with semi-hard triplet loss. Bottom Right: CNN with centroid loss.}
\label{fig:ContrastiveLearningInvariants}
\end{figure}

We train four separate student networks for the embeddings produced by the MLP and CNN teacher networks when trained with semi-hard triplet loss and centroid loss. We choose small students, consisting of two layers with 64 neurons and tanh activation.\footnote{For the MLP embeddings, the students are essentially as large as the teacher, which was a tiny NN; for the CNNs, the teacher was two orders of magnitude bigger than the students.} To check robustness, we trained the same student network architecture multiple times. We also changed the student architecture and activation functions, and generated two different embeddings by training different teacher networks. Our results are robust against these changes, and we present our findings in Figure~\ref{fig:ContrastiveLearningInvariants}.

We observe that several of the known knot invariants we tested seem to be encoded in what the teacher NNs learned, albeit to varying degree. As a baseline, we check how well the students can learn the teacher's learned map from the braid words to the embedding space. We see that the students can approximate this map with almost zero error. For the knot invariants constructed from the braid words, we find that the Goeritz matrix seems to be the simplest one to reproduce the teacher embeddings. The Alexander polynomial seems to not be part of the teacher encodings. This is perhaps a bit surprising, since the Alexander Polynomial can be constructed from Knot Floer Homology, which the teacher seems to learn to some extent. Moreover, the Alexander Polynomial (evaluated at -1) is given by the determinant of the (reduced) Goeritz matrix.

\subsubsection{Extracting invariants from the Transformer}
\begin{figure}[t]
\centering
\includegraphics[width=.6\textwidth]{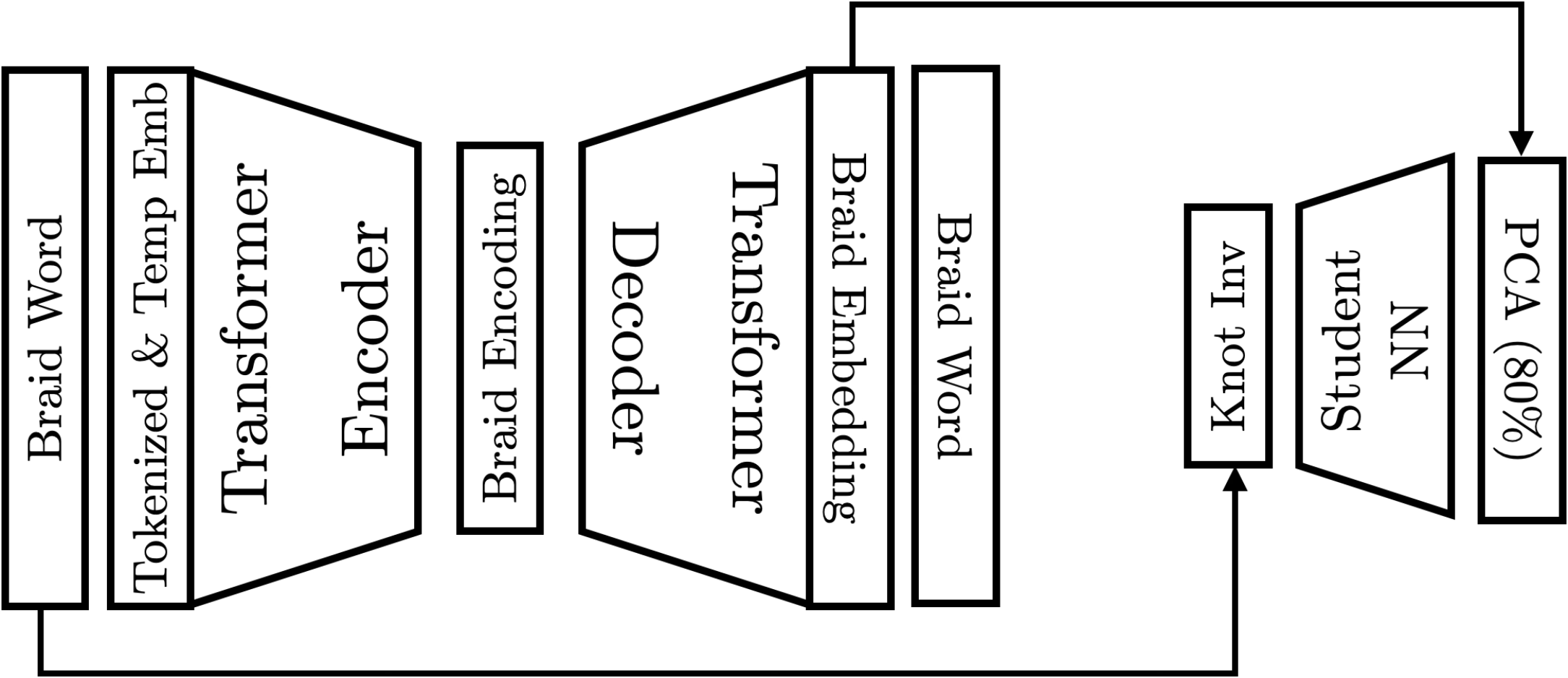}
\caption{Teacher-Student setup for the Transformer to extract which knot invariants are learned.}
\end{figure}

We use the output of the embedding layer, concatenate the embedding vectors of the output tokens, compute the principal components, and retain the first 200, which amounts to about 80 percent explained variance, and use these as labels for the student. The inputs (knot invariants) to the student are computed from braid words corresponding to the embeddings. To check how robust the invariants learned by the teacher are, we generate three sets of labels from three different teachers. They are trained independently with different input braids of size 25 and 30 for different numbers of epochs. The training time for the teachers varies between 2 to 4 hours on a NVIDIA GeForce RTX3090.

For the students, we can test robustness more extensively, since they are much cheaper to train. We focus on student NNs that are about a quarter the size of the teacher, so around 2.5M parameters. We tried different feed-forward NNs with some hyperparameter tuning, but found that the qualitative results are rather robust against network architecture and hyperparameter choices. We present the results for a four-layer NN with width 1024, batch normalization after each layer, and GELU activation function. We train them for 800 epochs with an Adam optimizer with pytorch~\cite{NEURIPS2019_9015} standard parameters on the various knot invariants as inputs, and the first 200 principal components of the embedding layer of the teacher as labels, as explained in Section~\ref{sec:idea_invariants}.

\begin{figure}[t]
\centering
\includegraphics[width=\textwidth]{./figs/LegendRuns}\\[-.2cm]
\subfloat[Teacher Run 1.]{\includegraphics[width=.32\textwidth]{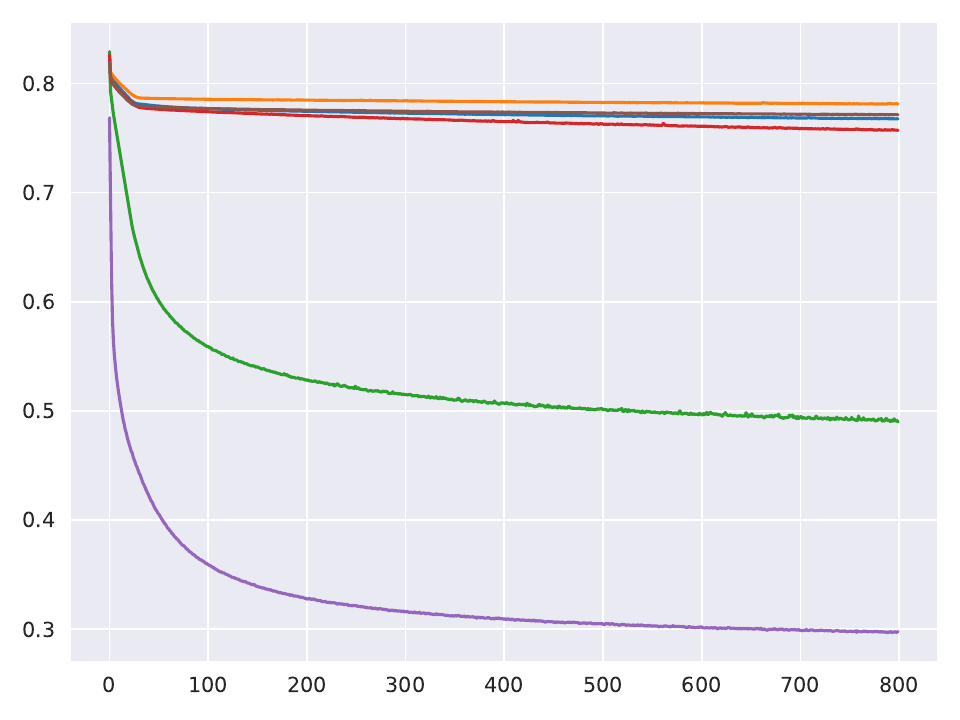}}~
\subfloat[Teacher Run 2.]{\includegraphics[width=.32\textwidth]{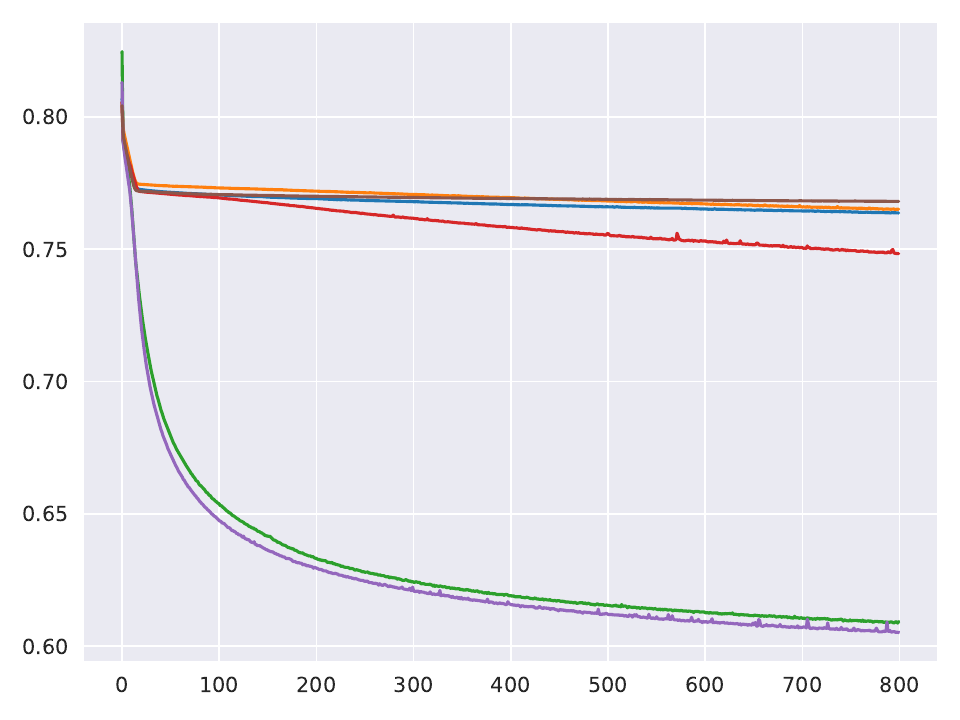}}~
\subfloat[Teacher Run 3.]{\includegraphics[width=.32\textwidth]{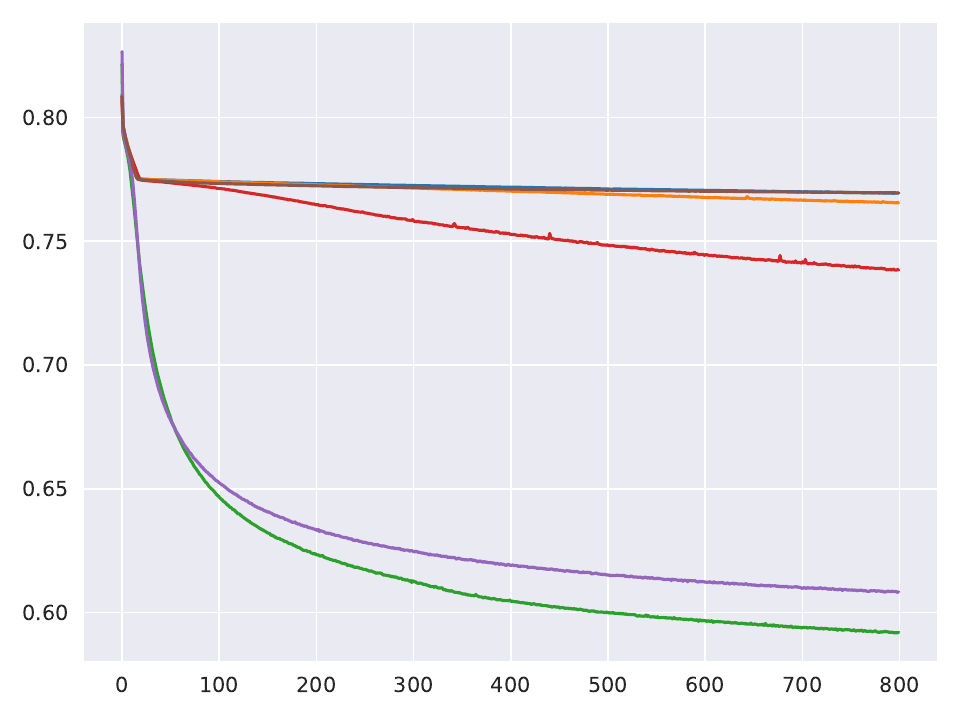}}
\caption{Loss of student network learning to predict different teacher embeddings from knot invariants.}
\label{fig:StudentTeacherInvariants}
\end{figure}

The results of the runs are given in Figure~\ref{fig:StudentTeacherInvariants}. We see that in each case, the student learns again to predict the teacher's embedding from the braid word and the Goeritz Matrix and, to a much lesser extent, from the homological data. This result is robust across all different student and teacher networks, as well as different distributions of input knots that we examined. In these experiments, the loss for the braid words should serve as a baseline against which one should compare the other losses; it is clear that the braid word can be mapped to the embeddings, since this is what the teacher learned. However, we should bear in mind that the loss on the braid words is not expected to be zero, since we are only using the first 200 PCs rather than the full embedding. Also, in contrast to the small, simple teachers in the contrastive learning setup, the Transformers are likely to learn more complex functions.

It is interesting that the students can learn to predict the embedding from the Goeritz matrix (almost) as well as from the braid word, which indicates that the Transformer is learning internally an embedding based on the Goeritz matrix without using other knot invariants that humans have developed. In particular, it is interesting that the teacher NN is not using the knot polynomials at all. For the Alexander polynomial, this is again surprising.

\subsection{Jones Unknot Conjecture\label{sec:exp_jones_unknot}}
In this section, we present the results of our search for knots that are not the unknot but nevertheless have Jones polynomial $J(K)=1$. We generated 400k knots with the sampling methods outlined in Section~\ref{sec:idea_jones_unknot}, but did not find a counterexample.

First, we generated a training set of $10^6$ inequivalent braid words of length 30 and computed their corresponding Jones polynomials using SageMath. To construct the braid words, we used the braid generator of~\cite{gukov2020learningunknot}, which essentially draws integers $\sigma_i\in[-N, N]$ from a random uniform distribution, resulting in a braid with up to $N+1$ strands (we choose $N=4$). We then ensure that the result is a one-component link. In addition, we make sure that a 30-crossing representative of the five smallest knots $0_1, 3_1, 4_1, 5_1, 5_2$ is in the training set. We had to add them by hand and use randomized braid identities and Markov moves to obtain a 30-crossing representation with 4 strands. We also made sure that the labels are unique, meaning that we did not generate two knots with the same Jones polynomial. We extracted the Jones polynomial coefficients, with the lowest term appearing being $a_{-78}t^{-78}$ and the highest term being $a_{76}t^{76}$, and represented them as a zero-padded feature vector of length $155$.

We then trained an encoder-decoder Transformer using the Int2Int package. We used a 256-dimensional embedding for the encoder tokens, and 4 layers with 8 attention heads for the encoder and decoder, which results in encoder and decoder models with about 5M parameters each. The dictionary size is around 1000, since we encode integers with base 1000. We trained the model for 450 epochs (where an epoch is defined as 300k training samples) with batch size 32, which takes about 3 days on a GeForce RTX4090.

With this trained model, we generated 100k knots with each of the sampling methods of Section~\ref{sec:idea_jones_unknot}. As a baseline against which we can compare the distribution of the knots generated from the decoder, we also generated 100k random braid words as outlined above. For the noise and temperature sampling, we find that if we increase the noise or temperature too much, the decoder starts returning nonsensical braid words such as 
\begin{align}
w=[1,2,-98111]\,.
\end{align}
While this is a well-defined link, it would have about 98k unlinked, unknotted components. If we encounter such a situation, we replace it by $[1,2,-4]$, which is a two-component link (that simplifies to two unlinked unknots). In cases where the sampling results in a braid generator that is outside of the training set range $[-4,4]$, we observe that the next step in the autoregressive generation very often is the `end of sequence' token. For that reason, we find that for larger noise, the lengths of the generated braid words drop.

\begin{figure}[t]
\centering
\subfloat[$\sigma$ vs number of crossings.\label{fig:GaussianSampling}]{\includegraphics[height=.35\textwidth]{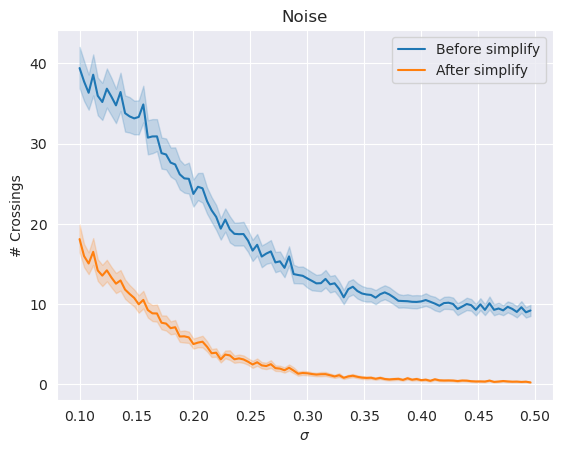}}\qquad\qquad
\subfloat[$T$ vs number of crossings.\label{fig:TempSampling}]{\includegraphics[height=.35\textwidth]{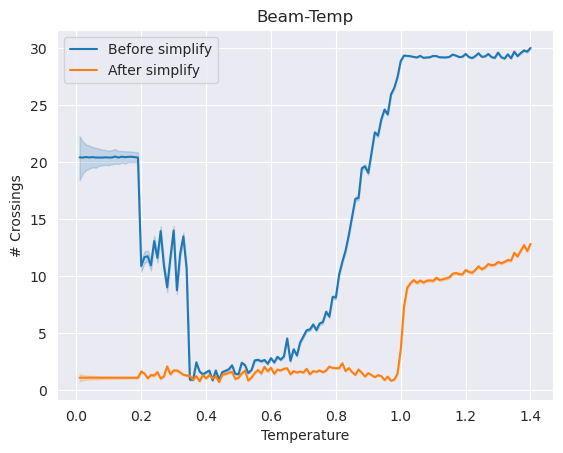}}
\caption{Dependence of the generated (blue) and generated+simplified (orange) braid words on the sampling parameter. We average over all braid words generated for this particular parameter value. Shaded regions are one standard deviation bands.}
\end{figure}

For Gaussian noise sampling, we increased the noise $\sigma$ linearly every 1000 knots. We generate 100k knots with noise in the interval $\sigma\in[0.1,0.5]$. As explained above, the higher the noise, the more likely we get an out-of-distribution prediction and the braid word is terminated subsequently, leading to shorter average lengths of the braids with increasing temperature. Gaussian noise with $\sigma<0.1$ predominantly produces unknots, which is why we chose this value as a lower bound. We observe that the average length of the generated braid words is around 40 for small $\sigma$, which is actually larger than the length of braid words in the training set, which is 30. Also, the generated braid words can be simplified to have on average only half as many crossings as the generated braid word. This could be a remnant of sampling around the unknot represented as a 30-crossing knot, which also allows for a lot of simplifications. To illustrate what kind of braid words we get, we present five randomly selected ones that were generated with $\sigma=.5$:
{\small
\begin{align*}
\begin{split}
w=&(8, -14, -7, -16, -5, -3, -2, 3, -9, -3, 18, -8, 4, 1, 10, -19, -13, 6, -12, -20, -17, -2, 15, 11)\\
w=&(-4, 10, -8, -2, -5, 1, -7, 9, -3, 2, 3, 6)\\
w=&(2, -3, 1, 4, 5)\\
w=&(-15, -26, -19, -17, -13, -8, -11, -10, 25, -22, 16, -21, 23, 9, -18, -7, -4, -24,\\& -5, 28, -14, 27, 2, -20, -1, 3, 6, 12)\\
w=&(-1, -3, 2, -5, -4, 6)
\end{split}
\end{align*}
}\noindent 
We also plot a few in appendix~\ref{app:GeneratedLinks}. In Figure~\ref{fig:GaussianSampling}, we show the average braid length for the generated (blue) and simplified (orange) braid words as a function of the variance $\sigma$.

Let us discuss path sampling next. Since we have a 256-dimensional embedding space, there are 256 trajectories $\gamma_1$ and 256 trajectories $\gamma_2$. We sample 10 points on each of the 256 trajectories, with 3 points along the $(t,0,\ldots,0)$ direction and 7 points along the $(0,t,\ldots,t)$ direction, for a combined total of $5120$ braid words. We then sample 1000 trajectories (out of 32640 possibilities) of type $\gamma_3$ and 1000 trajectories of type $\gamma_4$ by randomly choosing positions $(a,b)$. We again sample 10 times along these trajectories, with 3 points along the directions $(t,t,0,\ldots,0)$ and  7 points along the directions $(0,0,t\ldots,t)$, for a total of another $20k$ braids. We repeat this for trajectories from the unknot to the four knots with the lowest crossing number, resulting in a total of about $100k$ braids. The braid words generated along a path from the unknot to the trefoil are
{\small
\begin{align*}
\begin{split}
w=&(-2, -1, -2, 1, 2, -1, -3, -2, -3, -1, 2, 1, 2, 3, -3, -2, -1, -2, 1, 3, 2, 3, 1, -2, -1, 2, 1, 2, 4, -4)\\
w=&(-2, -1, -2, 3, 3, 2, 3, -2, -3, -3, -2, 1, 2, -1, -2, 1, 2, 3, 3, -2, -3, -4, -4, -4, 2, 1, 2, 3, 4, -3)\\
w=&(-3, -2, -1, -3, 2, 2, 2, -3, -2, -1, 2, 3, 1, 2, -3, -2, -1, -2, -2, 1, 2, 3, -2, -3, 2, 4, 4, 3, -4, 2)\\
w=&(-3, -3, -2, -1, 2, 3, 1, 1, 3, 2, 2, 2, -1, -2, -3, -1, -3, -2, 4, -3, -4, 3, 2, -3, -3, -3, 4, 2, 4, 4)\\
w=&(-3, -3, -3, -2, 3, 3, 2, -1, -1, -1, -2, 3, 1, 1, 2, 3, -4, -3, 4, -3, -3, 1, -2, -3, 4, 3, 4, 4, 2, -1)
\end{split}
\end{align*}
}\noindent 
Note that in contrast to the ones sampled with Gaussian noise, the results are much more in-distribution, in that they are all braid words of the same length with 5 strands. We again plot these and some more in Appendix~\ref{app:GeneratedLinks}.

For temperature plus beam sampling we linearly increase the temperature $T$ in the interval $T\in[0.01,1.40]$. At each temperature, we generate 100 braid words via a beam search, for a total of $100k$ braids. For temperatures below 0.01, we mostly get the unknot, which is why we start at $0.01$. We find that for temperatures up to $T=0.2$, the model mostly produces 20-crossing braid words which, however, can be simplified to knots with only a few crossings. In fact, most of the knots in the beam are equivalent and differ only by Markov moves and braid relations. There is a sharp transition around $T=0.2$ and $T=0.4$, where the average length of the generated braid words drops to 10 and 3, respectively. For even larger $T$, the braid word length goes up, which is trailed by the length of the simplified braids until $T=1$, where another phase transition seems to be happening. One explanation for the jumps at small $T$ is that the decoder keeps generating the same few knots in the vicinity of the unknot. Since the unknot in the training set is a complicated representation of a zero-crossing knot, the knots in the vicinity share this property, explaining the big discrepancy between the unsimplified and simplified braid words, cf.\ Figure~\ref{fig:TempSampling}. A few of the different braid words in a beam at temperature $T=.2$ look like this:
{\small
\begin{align*}
\begin{split}
w=&()\\
w=&(1, -2, -3, 2, 1, 3, 2, 3, -2, -3, 4)\\
w=&(-2, -1, -2, 1, 2, -1, -3, -2, -3, -1, 2, 1, 2, 3, -3, -2, 1, 1, 4)\\
w=&(-2, -1, -2, 1, 2, -1, -3, -2, -3, -1, 2, 1, 2, 3, -3, -2, -1, -2, 1, 3, 2, 3, 1, -2, -1, 2, 1, 2, 4, -4)\\
w=&(-2, -1, -2, 1, 2, -1, -3, -2, -3, -1, 2, 1, 2, 3, 3, -2, -1, -2, 1, 3, 2, 3, 1, -2, -1, 2, 4, 4)
\end{split}
\end{align*}
}\noindent 
At this temperature, they are also more in-distribution as the braid words generated via random Gaussian noise in the embedding space. As before, knot diagrams can be found in Appendix~\ref{app:GeneratedLinks}.

\begin{figure}[t]
\centering
\includegraphics[width=.9\textwidth]{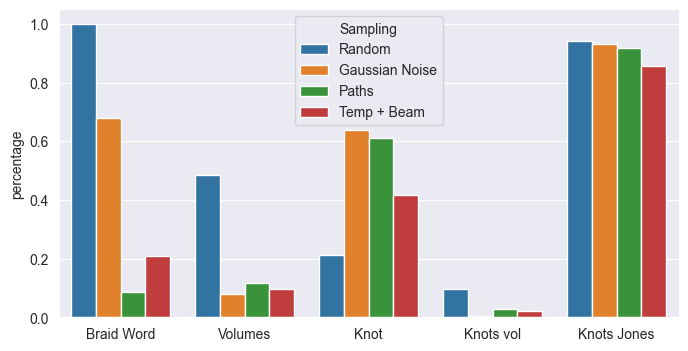}
\caption{Percentage of knots or links that are unique with respect to a given statistic, for different knot / link generation methods.}
\label{fig:SummaryGeneration}
\end{figure}

We compare the distribution of knots generated with these different sampling methods in Figure~\ref{fig:SummaryGeneration}, using different summary statistics over all 100k braid words generated by each method. First, we look at how many braids are different at the level of the braid word. For randomly generated braid words, all are different. For the other methods, we see that the number of different braid words is much lower; it is the lowest for path sampling, indicating that the knots between which we interpolate are indeed close together in the embedding space such that we do not encounter too many different braids along the trajectories. The fact that there are only very few inequivalent braid words for temperature + beam search as compared to Gaussian Noise sampling shows that across large temperature ranges, the braid words in the beam do not vary too much.

Next, we compute the volumes as a proxy for the number of inequivalent links in the sample. We see that all three different generation method result in roughly the same percentage of inequivalent hyperbolic links. However, which cannot be seen from the plot, is that the number of hyperbolic links (the volume is defined for these type of links) is lower for the generated braid words as compared to randomly sampled ones. For the latter, essentially all knots are hyperbolic, while the percentage is around 80\% for the generated ones. Out of these 80\% for which the volume is defined, around 90\% braid words turn out to describe the same knot or link. 

Next, we check how many of the links are actually knots, meaning they have only one component. For this statistic, we remove unknotted unlink components from braids. We find that their number is much higher as compared to the randomly generated ones. We also repeat the volume analysis for the one-component links only. For these, almost all generated knots turn out to be equivalent. Lastly, we compare the number of inequivalent Jones polynomials for the set of one-component knots that have different volumes. We see that all methods result in about 10\% of the knots having the same Jones polynomial.

\section{Conclusions}
\label{sec:conclusions}
We develop a framework for answering the central question in topology: when are two objects related by a geometric deformation? One goal of this study is to automate the task, so that neither theorems about equivalences nor knowledge of topological invariants are required. We train only on the fundamental data, i.e., different topological classes with multiple equivalent representatives in each class.

We illustrate our framework in the context of knot theory, where the central question is whether two knots are equivalent up to ambient space isotopy. We choose braid words to represent the data, since they have been found to outperform other representation for machine learning of topological invariants~\cite{lindsay2025learnabilityknotinvariantsrepresentation}. Moreover, braid words bear resemblance to natural language and hence allow us to use recent advances in natural language processing to generate new knots in the same class.

Using different losses in a contrastive learning setup as well as a generative setup, we train neural networks to learn topology-invariant embeddings of their input data: different representatives of the same equivalence class cluster in the embedding space, and clusters of knots with different topology are separated. In this way, the network is trained to ``abstract away'' information in the input data that is not pertaining to topology, and the position of the cluster in the embedding space gives a machine-learned topological invariant that need not have strong correlations with known topological invariants. For the contrastive learning setup, we use semi-hard triplet loss, and we develop a new loss based on centroids of clusters. We employ CNNs that make the network invariant with respect to one of the symmetries present in braid representations of knots.  For the generative task, we train an encoder-decoder Transformer pair, where we have different classes of knots as inputs, and each class has multiple representatives. For each class, we choose a random representative as a label and train the generator to map all inputs from the same topology class to this representative. To measure accuracy of the generator, we identify the topology class of the generated knot by comparing the hyperbolic knot volume.

We employ a student-teacher setup to interpret which knot invariants were learned during the training process. The student network gets known knot invariants as inputs and is tasked to learn a map to the embedding space of the teacher. If a small student NN can learn this map to high accuracy, we postulate that the corresponding knot invariant is (partially) learned by the teacher in computing the embedding. Applying this procedure to the embeddings learned by teacher NNs trained with contrastive learning as well as in the generative setup, we found that all networks seem to utilize the Goeritz matrix and to a much lesser degree homological data (Knot Floer Homology). The NNs trained with contrastive learning also learned scalar invariants like the hyperbolic knot volume, the $\nu$- and $\tau$-invariant, and the Jones Polynomial. None of the teacher NNs seem to make use of the Alexander Polynomial.

Lastly, we used the auto-regressive encoder-decoder Transformers to generate potential counterexamples to the Jones Unknot Conjecture. We train the NN to learn a map from the Jones Polynomial to a braid representation of a knot with this Jones Polynomial. We then search the embedding space of the encoder around the unknot, sample trajectories from the unknot to other knots with simple Jones Polynomials, and perform a beam search paired with increasing the temperature at the output layer to produce knots that are distinct from the unknot but still have trivial Jones polynomials. We generate several 100k knots, but none of them have trivial Jones Polynomial.

There are a number of natural directions for future work. First, one can combine the contrastive and generative methods into one to facilitate training or precondition the embedding on certain topological invariants. This could improve performance in a supervised learning setup where networks should learn to predict topological invariants, or in the generation of potential counterexamples with certain properties. Second, while we  have applied our framework in the context of knot theory, it is much more general and could be applied to other problems, such as the classification of three-manifolds. Third, as explained in Section~\ref{sec:idea_invariants}, one can subtract the embeddings learned by the student from those learned by the teacher such that the resulting embedding pertains to new topological invariants that are learned by the teacher and do not correspond to any known invariants. One can then apply symbolic regression or FunSearch~\cite{Romera-Paredes:2024aaa} to extract how these invariants are computed. Finally, synthetic data, such as those we have studied, are a useful avenue for better understanding ML theory. Specifically, the dynamics of learned features in the embedding dimension and how they are related to features of the data is an important open problem. In many ML applications, the features of data are ill-defined from a mathematical perspective; for instance, handwritten digits are classified according to how humans would label them, and may be misclassified when a digit appears to be a blend of two clear digits. Synthetic data avoids this uncertainty and accordingly the rigorously defined features of the data in principle allow for better comparison to learned features in the embedding dimension.

\vspace{1cm}
\noindent \textbf{Acknowledgements.} We are grateful to Fran\c{c}ois Charton, Alex Davies, Mike Douglas, Mike Freedman, Sergei Gukov, Mark Hughes, and Onkar Singh Gujral for discussions regarding this project. J.H.\ is supported by NSF CAREER grant PHY-1848089. F.R.\ is supported by NSF grant PHY-2210333 and startup funding from Northeastern University. Both authors are supported by the National Science Foundation under Cooperative Agreement PHY-2019786 (The NSF AI Institute for Artificial Intelligence and Fundamental Interactions).

\appendix
\clearpage
\section{Appendix: Generated links}
\label{app:GeneratedLinks}
\begin{figure}[h!t]
\centering
\includegraphics[width=\textwidth]{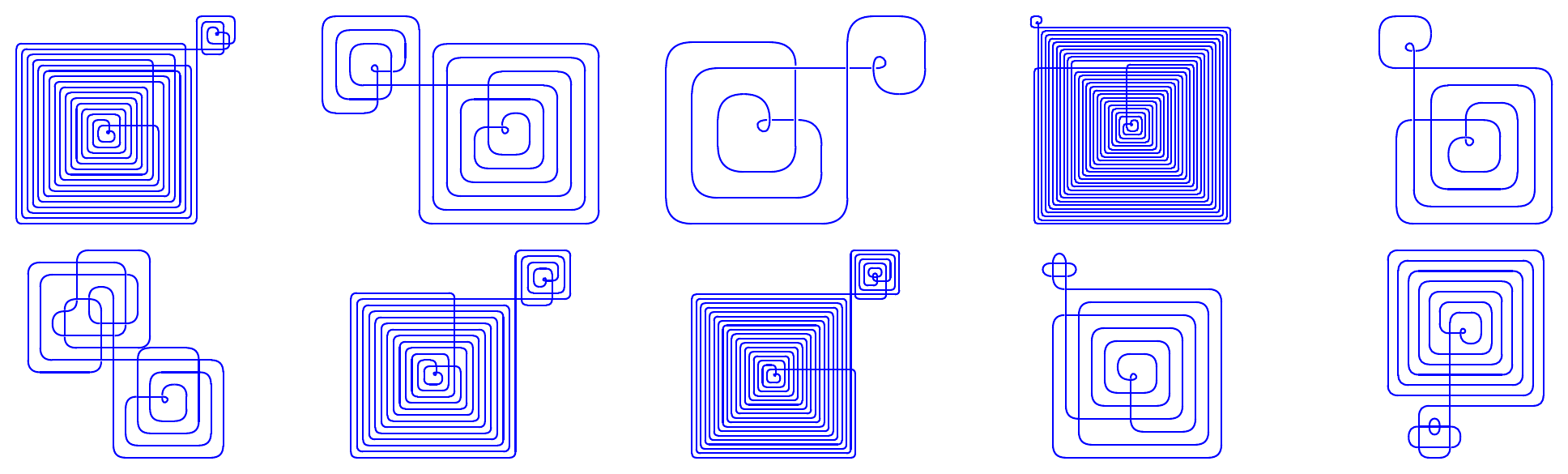}
\rule{\textwidth}{0.4pt}\\[12pt]\noindent
\includegraphics[width=\textwidth]{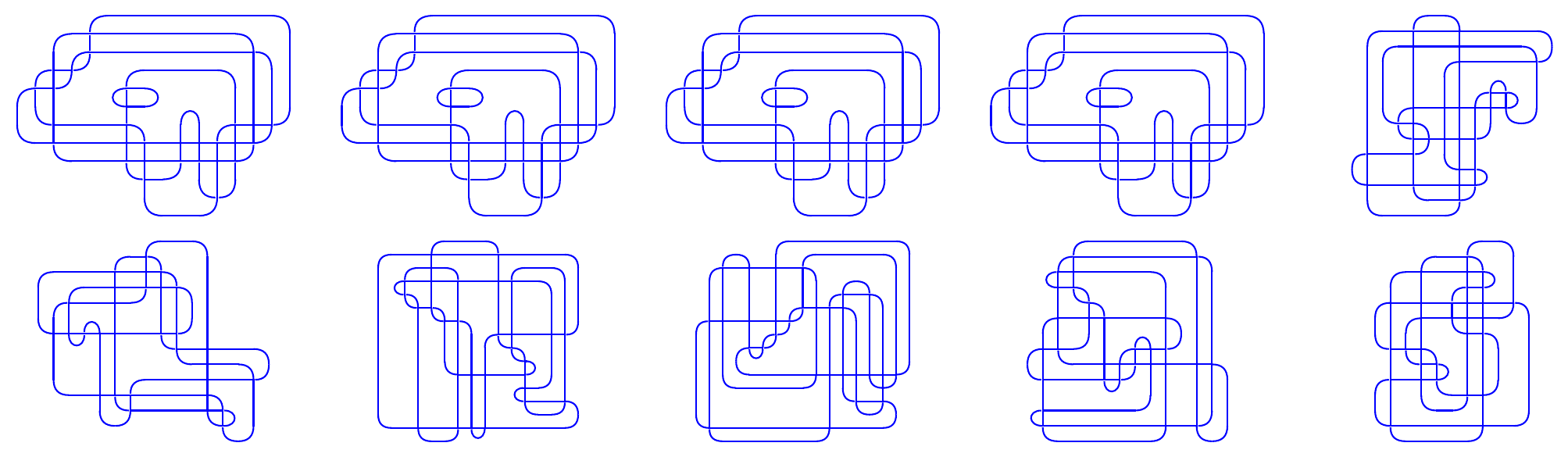}
\rule{\textwidth}{0.4pt}\\[12pt]\noindent
\medskip
\includegraphics[width=\textwidth]{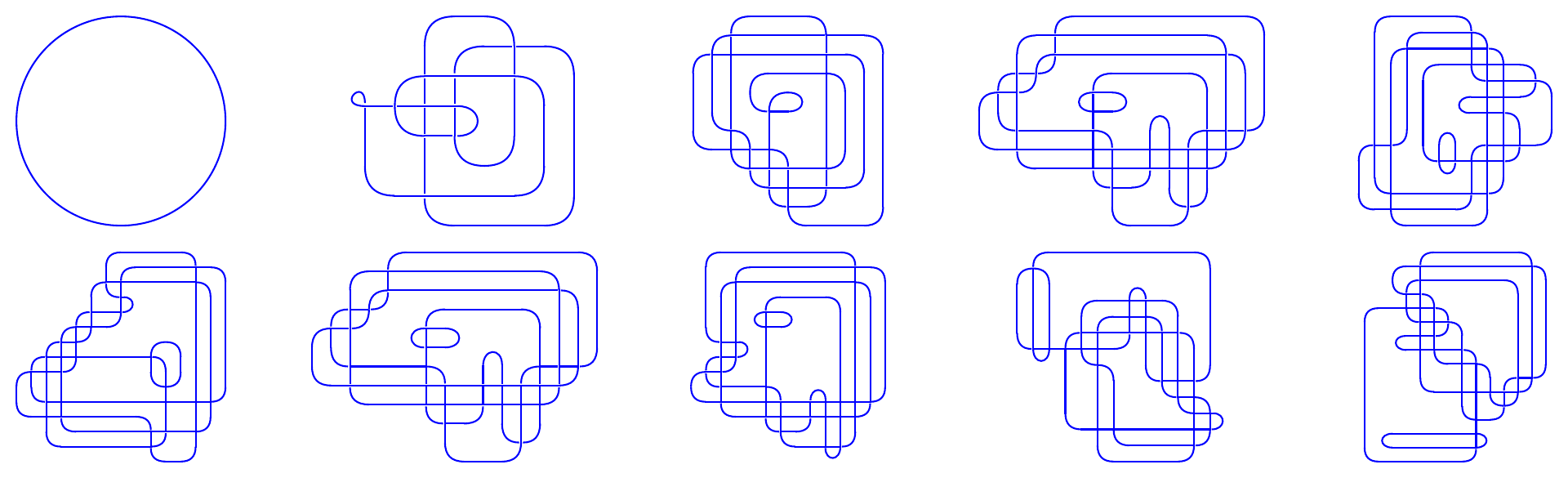}
\caption{Top: 10 knots generated with Gaussian noise ($\sigma=.2$). Center: 10 knots sampled along a trajectory from the Unknot to the Trefoil. Bottom: 10 knots in a beam search at $T=1.0$.}
\end{figure}
\clearpage

\begin{landscape}
\section{Appendix: Beam Search of Braids}\label{app:BeamSearch}
    \begin{table}[ht]
    \centering
    {
    \begin{tabular}{|c|c|c|}
    \hline
    \footnotesize\textbf{\text{vol} 33.69} & \footnotesize\textbf{Braid representation} & \footnotesize\textbf{Comment} \\
    \hline
    {\footnotesize Braid 1} & {\footnotesize $[-2, -3, -3, -3, -4, 1, -4, 3, -2, -2, -1, -1, -2, 4, -3, 2, -4, -3, 2, 2, 2, -1, 2, -3, -4, 3, -2, 1, 1, -3]$} &  {\footnotesize }\\
    \hline
    {\footnotesize Braid 2} & {\footnotesize $[\cb{-3}, \cb{-3}, {-4}, {-4}, 1, 3, -2, -2, -1, -1, -2, 4, -3, 2, -4, -3, 2, 2, 2, -1, 2, -3, -4, 3, -2, 1, 1, -3, \cb{-2}, \cb{-3}]$} &  {\footnotesize 3 $\times$ Conj. of Braid 1}\\
    \hline
    {\footnotesize Braid 3} & {\footnotesize $[-3, -3, -4, -4, 1, 3, -2, -2, -1, -1, -2, 4, -3, 2, -4, -3, 2, 2, 2, -1, 2, -3, -4, 3, -2, 1, 1, \cb{-2}, \cb{-3}, \cb{-2}]$} & {\footnotesize BR of Braid 2}\\
    \hline
    {\footnotesize Braid 4} & {\footnotesize $[\cb{-4}, 3, -2, -2, -1, -1, -2, 4, -3, 2, -4, -3, 2, 2, 2, -1, 2, -3, -4, 3, -2, 1, 1, \crr{-3}, \crr{-2}, \crr{-3}, \cb{-3}, \cb{-3}, \cb{1}, \cb{-4}]$} & {\footnotesize \crr{BR}, \cb{2 $\times$ Conj. + BR} of Braid 3} \\
    \hline
    \hline
    \footnotesize\textbf{\text{vol}: 38.2} & \footnotesize\textbf{Braid representation} & \footnotesize\textbf{Comment} \\
    \hline
    {\footnotesize Braid 1} & {\footnotesize $[1, -2, 3, 4, -3, 2, -3, -3, 1, -2, -4, 3, -4, -4, 3, 4, 1, 1, -3, 2, -3, 4, 4, -3, -2, -3, -3, -3, 1, 4]$} &  {\footnotesize }\\
    \hline
    {\footnotesize Braid 2} & {\footnotesize $[\cb{1}, {1}, -2, 3, 4, -3, 2, -3, -3, 1, -2, -4, 3, -4, -4, 3, 4, 1, 1, -3, 2, -3, 4, 4, -3, -2, -3, -3, -3, \cb{4}]$} & {\footnotesize Conj. + BR of Braid 1}\\
    \hline
    {\footnotesize Braid 3} & {\footnotesize $[{-3}, {-3}, {1}, {1}, -2, {4}, 3, 4, -3, 2, -3, -3, 1, -2, -4, 3, -4, -4, 3, 4, 1, 1, -3, 2, -3, 4, 4, -3, -2, -3]$} & {\footnotesize }\\
    \hline
{\footnotesize Braid 4} & {\footnotesize $[\cb{-3}, {-3}, {-3}, {1}, 1, -2, {4}, 3, 4, -3, 2, -3, -3, 1, -2, -4, 3, -4, -4, 3, 4, 1, 1, -3, 2, -3, 4, 4, -3, \cb{-2}]$} & {\footnotesize Conj. of Braid 3 }\\
    \hline
    \hline
    \footnotesize\textbf{\text{vol}: 31.2} & \footnotesize\textbf{Braid representation} & \footnotesize\textbf{Comment} \\
    \hline
    {\footnotesize Braid 1} & {\footnotesize $[1, -4, 2, 3, -2, 3, -4, -4, 1, 2, -3, 4, 4, 3, 3, -2, 3, -1, 4, 3, -2, -1, 3, 4, -2, 1, 4, 3, 2, -4]$} &  {\footnotesize }\\
    \hline
    {\footnotesize Braid 2} & {\footnotesize $[{3}, {-2}, {3}, {-1}, {4}, {3}, {-2}, 3, {-1}, {-2}, {4}, {1}, {4}, {3}, {-4}, {2}, {1}, {-4}, 2, 3, -2, 3, -4, -4, 1, 2, -3, 4, 4, 3]$} & {\footnotesize }\\
    \hline
    {\footnotesize Braid 3} & {\footnotesize $[3, -2, 3, -1, 4, 3, -2, 3, -1, -2, 4, 1, 4, 3, \cb{2}, \cb{1}, \cb{-4}, \cb{-4}, 2, 3, -2, 3, -4, -4, 1, 2, -3, 4, 4, 3]$} & {\footnotesize 4 $\times$ BR of Braid 2}\\
    \hline
    {\footnotesize Braid 4} & {\footnotesize $[1, -4, 2, 3, -2, 3, -4, -4, 1, 2, -3, 4, 4, 3, 3, -2, 3, -1, 4, 3, -2, -1, 3, 4, -2, 1, 4, 3, \cb{-4}, \cb{2}]$} & {\footnotesize BR of Braid 1}\\
    \hline
    \hline
    \footnotesize\textbf{\text{vol}: 28.4} & \footnotesize\textbf{Braid representation} & \footnotesize\textbf{Comment} \\
    \hline
    {\footnotesize Braid 1} & {\footnotesize $[3, 3, 2, -1, -2, 1, 1, 2, 1, -2, 1, 1, -3, -2, 1, -2, 4, -1, 2, 3, 1, -2, 3, 3, 1, 3, 3, 4, 3, 2]$} &  {\footnotesize }\\
    \hline
    {\footnotesize Braid 2} & {\footnotesize $[\cb{3}, {2}, {-1}, {-2}, 1, 1, 2, 1, -2, 1, 1, -3, -2, 1, -2, 4, -1, 2, 3, 1, -2, 3, 3, 1, 3, {3}, 4, 3, 2, \cb{3}]$} & {\footnotesize Conj. of Braid 1}\\
    \hline
    {\footnotesize Braid 3} & {\footnotesize $[3, 2, -1, -2, 1, 1, 2, 1, -2, 1, 1, -3, -2, 1, -2, 4, -1, 2, 3, 1, -2, 3, 3, 1, 3, \cb{4}, \cb{3}, \cb{4}, {2}, 3]$} & {\footnotesize BR of Braid 2}\\
    \hline
    {\footnotesize Braid 4} & {\footnotesize $[\cb{3}, 3, 2, -1, -2, 1, 1, 2, 1, -2, 1, 1, -3, -2, 1, -2, 4, -1, 2, 3, 1, -2, 3, 3, 1, 3, {4}, {3}, {4}, \cb{2}]$} & {\footnotesize Conj. of Braid 3}\\
    \hline
    \end{tabular}}
    \caption{Beam search results with noted differences for four different knot classes. Comments indicate the use of conjugation moves or braid relations on the stated braid, highlighting where the braid word differs from the reference braid.}
    \label{tab:beam_search_results}
    \end{table}
    
    \end{landscape}
\clearpage


\providecommand{\href}[2]{#2}

\end{document}